\documentclass[12pt]{amsart}
\numberwithin{equation}{section}

\usepackage{accents}
\usepackage{appendix}
\usepackage{amsfonts}
\usepackage{amsmath}
\usepackage{amssymb}	
\usepackage{amsthm}
\usepackage{array,booktabs,multirow} 
\usepackage{braket}
\usepackage{cite}
\usepackage{dsfont}
\usepackage{graphicx}
\usepackage{mathalfa}
\usepackage{mathtools}
\usepackage[shortlabels]{enumitem} 
\usepackage{etoolbox}
\usepackage{float}
\usepackage[hang, flushmargin]{footmisc}
\usepackage{latexsym}
\usepackage{lipsum}
\usepackage{needspace}
\usepackage{tikz}
\usepackage{hyperref}
\usetikzlibrary{matrix,arrows}  
\usepackage{multicol}

\newcommand{\dbar}{\,\mathchar'26\mkern-12mu d \hspace{0.06em}}

\theoremstyle{plain}
\newtheorem{thm}{Theorem}[section]

\newtheorem{lem}[thm]{Lemma}
\newtheorem{prop}[thm]{Proposition}

\theoremstyle{definition}

\newtheorem{defn}[thm]{Definition}

\newtheorem{prob}[thm]{Problem}

\theoremstyle{remark}

\AtBeginEnvironment{lem}{\Needspace{3\baselineskip}}%
\AtBeginEnvironment{thm}{\Needspace{3\baselineskip}}
\AtBeginEnvironment{prop}{\Needspace{3\baselineskip}}
\AtBeginEnvironment{cor}{\Needspace{3\baselineskip}}

\setlist[enumerate,1]{leftmargin=2em}

\def\F{\mathbb F}
\def\I{\mathbf I}
\def\P{\mathbf P}
\def\S{\mathfrak S}
\def\N{\mathbb N}

\def\e{\varepsilon}
\def\z2s4{{}_{\,\{\pm 1\}}\!\!\simeq_{\,\S_4}}
\def\fz2s4{{}_{\,\{\pm 1\}}\!\!\approx_{\,\S_4}}

\title[Verma $\triangle_q$-modules and irreducible $\triangle_q$-modules at roots of unity]{Verma modules and finite-dimensional irreducible modules of the universal Askey--Wilson algebra \\
at roots of unity}

\author{Hau-Wen Huang}
\address[H.-W. Huang]{
Department of Mathematics\\
National Central University\\
Chung-Li 32001 Taiwan
}
\email{hauwenh@math.ncu.edu.tw}

\begin{document}

\begin{abstract}
Assume that $\mathbb F$ is an algebraically closed field and fix a nonzero scalar $q\in \mathbb F$ with $q^4\not=1$. The universal Askey--Wilson algebra $\triangle_q$ is a unital associative algebra over $\mathbb F$ defined by generators and relations. The generators are $A,B,C$ and the relations assert that each of 
\begin{gather*}
A+\frac{qBC-q^{-1}CB}{q^2-q^{-2}},
\qquad 
B+\frac{qCA-q^{-1}AC}{q^2-q^{-2}},
\qquad 
C+\frac{qAB-q^{-1}BA}{q^2-q^{-2}}
\end{gather*}
commutes with $A,B,C$. The Verma $\triangle_q$-modules are a family of infinite-dimensional $\triangle_q$-modules with marginal weights. 
Under the condition that $q$ is not a root of unity, 
it was shown that every finite-dimensional irreducible $\triangle_q$-module has a marginal weight and is isomorphic to a quotient of a Verma $\triangle_q$-module. 
Assume that $q$ is a root of unity. We prove that every finite-dimensional irreducible $\triangle_q$-module with a marginal weight is isomorphic to a quotient of a Verma $\triangle_q$-module. Properly speaking, two natural families of finite-dimensional quotients of Verma $\triangle_q$-modules contain all finite-dimensional irreducible $\triangle_q$-modules with marginal weights up to isomorphism. 
Furthermore, we classify the finite-dimensional irreducible $\triangle_q$-modules with marginal weights up to isomorphism.
\end{abstract}

\maketitle



\section{Introduction}\label{s:intro}
Throughout this paper, we adopt the following conventions:
A vacuous product is equal to the multiplicative identity. A vacuous summation is equal to the additive identity. Let $\N$ denote the set of all nonnegative integers.
Assume that the underlying field $\F$ is algebraically closed.  
Let $\F^\times$ denote the multiplicative group of all nonzero scalars in $\F$. 
Fix a scalar $q\in \F^\times$ with $q^4\not=1$. 
Let $x$ denote an indeterminate over $\F$. 
For any $a\in \F$ let $\sqrt{a}$ denote a fixed root of $x^2-a$ in $\F$. 
Given a left (resp. right) action of a group $G$ on a set $S$, 
the notation $G\backslash S$ (resp. $S/G$) stands for the set of all $G$-orbits of $S$.

The Askey--Wilson algebras \cite{hidden_sym,LP&dual} are 
a family of unital associative algebras defined by generators and relations. 
These algebras describe the bispectral properties of orthogonal polynomials in the Askey scheme \cite{Koe2010}.  
Since the advent of Askey--Wilson algebras, they have been found to have applications in various fields such as $P$- and $Q$-polynomial association schemes \cite{BannaiIto2021,Lthypercube,qRacahDRG,hypercube2002,halved:2023,halved:2024,Huang:CG&Hamming,Huang:CG&Johnson,Huang:CG&Grassmann,odd:2024}, spin models \cite{cur2007-1,nomura2021,cur2007-2,nomura2025}, 
Leonard pairs \cite{Askeyscheme,lp&awrelation,Vidunas:2007,Huang:2012},  
Lie algebras \cite{uaw&equit2011,lp&equ2011,
CG2013,gvz2013,SH:2020,cur2004-1,K&sl2}, 
double affine Hecke algebras \cite{Huang:AW&DAHAmodule,daha&Z3,koo07,koo08,daha&AW
,Huang:R<BImodules,Huang:R<BI,BI&NW2016,daha&LP}, coupling problems \cite{Levy1965,gz93,Huang:CG,Huang:RW} and superintegrable systems \cite{BI2014-2,R&LD2014,gvz2014}.
The universal Askey-Wilson algebra is a central extension of the Askey-Wilson algebras corresponding to the most general orthogonal polynomials in the Askey scheme, namely the Askey--Wilson polynomials and the $q$-Racah polynomials. The definition is presented as follows:

\begin{defn}
[Definition 1.2, \cite{uaw2011}]
The {\it universal Askey--Wilson algebra} $\triangle_q$ is a unital associative algebra over $\F$ defined by generators and relations. The generators are $A,B,C$ and the relations assert that each of 
\begin{gather}
\label{abc}
A+
\frac{qBC-q^{-1}CB}{q^2-q^{-2}},
\qquad
B+
\frac{qCA-q^{-1}AC}{q^2-q^{-2}},
\qquad
C+
\frac{qAB-q^{-1}BA}{q^2-q^{-2}}
\end{gather}
is central in $\triangle_q$. 
\end{defn}

Let $\alpha,\beta,\gamma$ denote the central elements of $\triangle_q$ obtained by multiplying the elements (\ref{abc}) by $q+q^{-1}$ respectively. In other words 
\begin{eqnarray}
\frac{\alpha}{q+q^{-1}}
&=&
A+\frac{qBC-q^{-1}CB}{q^2-q^{-2}},
\label{e:alpha}
\\
\frac{\beta}{q+q^{-1}}
&=&
B+
\frac{qCA-q^{-1}AC}{q^2-q^{-2}},
\label{e:beta}
\\
\frac{\gamma}{q+q^{-1}}
&=&
C+
\frac{qAB-q^{-1}BA}{q^2-q^{-2}}.
\label{e:gamma}
\end{eqnarray}

\begin{prop}
\label{prop:UAW_presentation}
The algebra $\triangle_q$ has a presentation given by generators $A,B,\alpha,\beta,\gamma$ and the relations assert that $\alpha,\beta,\gamma$ are central in $\triangle_q$ and 
\begin{gather}
\alpha = 
\frac{B^2A-(q^2+q^{-2})BAB+AB^2+(q^2-q^{-2})^2A+(q-q^{-1})^2B\gamma}
{(q-q^{-1})(q^2-q^{-2})}, \label{e:a->ABc}
\\
\beta =\frac{A^2B-(q^2+q^{-2})ABA+BA^2+(q^2-q^{-2})^2B+(q-q^{-1})^2A\gamma}{(q-q^{-1})(q^2-q^{-2})}. \label{e:b->ABc}
\end{gather}
\end{prop}
\begin{proof}
The relations (\ref{e:a->ABc}) and (\ref{e:b->ABc}) are obtained by applying (\ref{e:gamma}) to eliminate $C$ in (\ref{e:alpha}) and (\ref{e:beta}).
\end{proof}

Let $V$ denote a $\triangle_q$-module. 
For any $\mu\in \F^\times$ define
$$
V(\mu)=\{v\in V\,|\, (B -\mu-\mu^{-1})v=0\}.
$$
Note that $V(\mu)=V(\mu^{-1})$ for any $\mu\in \F^\times$. 
A scalar $\mu\in \F^\times$ is called a {\it weight} of $V$ whenever $V(\mu)\not=\{0\}$. In this case $V(\mu)$ is called a {\it weight space} of $V$ with weight $\mu$ and every nonzero $v\in V(\mu)$ is called a {\it weight vector} of $V$ with weight $\mu$.

%\begin{lem}
%\label{lem0:weight}
%For any $\triangle_q$-module $V$ and any weight $\mu$ of $V$, the element 
%$$
%(B-\mu-\mu^{-1})(B-\mu q^2-\mu^{-1} q^{-2})(B-\mu q^{-2}-\mu^{-1} q^2)A
%$$ 
%vanishes on $V(\mu)$.
%\end{lem}
%\begin{proof}
%Let $v\in V(\mu)$ be given.
%Applying $v$ to either side of (\ref{e:a->ABc}) 
%it follows that 
%\begin{gather}
%\label{BBAv}
%(B-\mu q^2-\mu^{-1} q^{-2})(B-\mu q^{-2}-\mu^{-1} q^2)A 
%v
%\end{gather}
%is equal to $(q-q^{-1})^2$ times 
%$
%(q+q^{-1})\alpha v-(\mu+\mu^{-1})\gamma v$. 
%Since $\alpha$ and $\gamma$ are central in $\triangle_q$ the vectors $\alpha v$ and $\gamma v$ are in $V(\mu)$. Therefore (\ref{BBAv}) lies in $V(\mu)$. The lemma follows.
%\end{proof}

\begin{lem}
\label{lem:weight}
For any $\triangle_q$-module $V$ and any weight $\mu$ of $V$ the following relations hold:
\begin{enumerate}
\item $(B-\mu q^2-\mu^{-1} q^{-2})(B-\mu q^{-2}-\mu^{-1} q^2)A V(\mu)\subseteq V(\mu)$.

\item $(B-\mu q^2-\mu^{-1} q^{-2}) (B-\mu-\mu^{-1}) 
AV(\mu)\subseteq V(\mu q^{-2})$.

\item $(B-\mu q^{-2}-\mu^{-1} q^{2}) (B-\mu-\mu^{-1}) 
AV(\mu)\subseteq V(\mu q^{2})$.
\end{enumerate}
\end{lem}
\begin{proof}
(i): Let $v\in V(\mu)$ be given.
Applying $v$ to either side of (\ref{e:a->ABc}) 
it follows that 
\begin{gather}
\label{BBAv}
(B-\mu q^2-\mu^{-1} q^{-2})(B-\mu q^{-2}-\mu^{-1} q^2)A 
v
\end{gather}
is equal to $(q-q^{-1})^2$ times 
$
(q+q^{-1})\alpha v-(\mu+\mu^{-1})\gamma v$. 
Since $\alpha$ and $\gamma$ are central in $\triangle_q$ the vectors $\alpha v$ and $\gamma v$ are in $V(\mu)$. Therefore (\ref{BBAv}) lies in $V(\mu)$. The relation (i) follows.

(ii), (iii): The relation (i) implies that 
$$
(B-\mu-\mu^{-1})(B-\mu q^2-\mu^{-1} q^{-2})(B-\mu q^{-2}-\mu^{-1} q^2)A V(\mu)=\{0\}.
$$ 
The relations (ii) and (iii) are immediate from the above equation.
\end{proof}

\begin{defn}
\label{defn:marginal}
A weight $\mu$ of the $\triangle_q$-module $V$ is said to be {\it marginal} 
if there exists a weight vector $v$ of $V$ with weight $\mu$ such that 
$$
(B-\mu q^2-\mu^{-1} q^{-2}) (B-\mu-\mu^{-1}) 
Av=0.
$$
\end{defn}

If $q$ is not a root of unity, it was shown that all finite-dimensional irreducible $\triangle_q$-modules have marginal weights and can be constructed from the following infinite-dimensional $\triangle_q$-modules up to isomorphism:

\begin{thm}
[Section 3, \cite{Huang:2015}]
\label{thm:Verma}
For any $(a,b,c,\lambda)\in {\F^\times}^4$ there exists an infinite-dimensional $\triangle_q$-module $M_\lambda(a,b,c)$ satisfying the following conditions:
\begin{enumerate}
\item There exists a basis $\{m_i\}_{i\in \N}$ for the $\triangle_q$-module $M_\lambda(a,b,c)$ such that 
\begin{align*}
(A-\theta_i) m_i&=m_{i+1} \qquad 
\hbox{for all $i\in \N$},
\\
(B-\theta_i^*) m_i&= \varphi_i m_{i-1}
\qquad 
\hbox{for all $i\in \N$},
\end{align*}
where $m_{-1}$ is interpreted as any vector of $M_\lambda(a,b,c)$ and 
\begin{align}
\theta_i &= a\lambda^{-1} q^{2i}+ a^{-1} \lambda q^{-2i}
\qquad 
\hbox{for all $i\in \N$},
\label{thetai}
\\
\theta_i^*&= b\lambda^{-1} q^{2i}+ b^{-1} \lambda q^{-2i}
\qquad 
\hbox{for all $i\in \N$},
\label{thetais}
\\
\begin{split}\label{varphii}
\varphi_i &=
a^{-1}b^{-1} \lambda q(q^i- q^{-i})(\lambda^{-1}q^{i-1}-\lambda q^{1-i})
\\
&\qquad \times \, (q^{-i}- abc\lambda^{-1}q^{i-1})(q^{-i} - abc^{-1}\lambda^{-1}q^{i-1})
\qquad 
\hbox{for all $i\in \N$}.
\end{split}
\end{align}

\item The elements $\alpha,\beta,\gamma$ act on $M_\lambda(a,b,c)$ as scalar multiplication by 
\begin{align}
%\omega&=
(b+b^{-1})(c+c^{-1})
+(a+a^{-1})(\lambda q+\lambda^{-1} q^{-1}),
\label{omega}
\\
%\omega^*&=
(c+c^{-1})(a+a^{-1})
+(b+b^{-1})(\lambda q+\lambda^{-1} q^{-1}),
\label{omega*}
\\
%\omega^\e&=
(a+a^{-1})(b+b^{-1})
+(c+c^{-1})(\lambda q+\lambda^{-1} q^{-1}),
\label{omegae}
\end{align}
 respectively.
\end{enumerate}
\end{thm}

By Theorem \ref{thm:Verma}(i) the $\triangle_q$-module $M_\lambda(a,b,c)$ has the marginal weight $b\lambda^{-1}$. 
In 2009 the present author thought up the rough idea for creating $M_\lambda(a,b,c)$ during his work \cite{Huang:2012} on Leonard triples of $q$-Racah type. Those Leonard triples provide a family of finite-dimensional irreducible $\triangle_q$-modules. 
In the 2015 paper \cite{Huang:2015}, the $\triangle_q$-module $M_\lambda(a,b,c)$ was formally launched to classify the finite-dimensional irreducible $\triangle_q$-modules at $q$ not a root of unity. The $\triangle_q$-module $M_\lambda(a,b,c)$ is called the {\it Verma $\triangle_q$-module} as the contribution of Verma modules did in the semisimple Lie algebras.

In this article, we focus on those finite-dimensional irreducible $\triangle_q$-modules with marginal weights under the condition that $q$ is a root of unity.
From now on, we always assume that $q$ is a root of unity  with order $d\not=1,2,4$ and set 
$$
\dbar=\left\{
\begin{array}{ll}
d  \qquad &\hbox{if $d$ is odd},
\\
\frac{d}{2} \qquad &\hbox{if $d$ is even}.
\end{array}
\right.
$$
Note that $\dbar$ is the order of $q^2$ and $\dbar\geq 3$. In \cite{Huang:2021} it was shown that the dimension of every finite-dimensional irreducible $\triangle_q$-module is less than or equal to $\dbar$. Moreover every irreducible $\triangle_q$-module of dimension less than $\dbar$ has a marginal weight.

There are two natural families of finite-dimensional $\triangle_q$-modules obtained by taking quotients of Verma $\triangle_q$-modules.
One of the two families was first released in \cite[Section 4]{Huang:2015} and the description is as follows: Pick a triple $(a,b,c)\in {\F^\times}^3$ and any $n\in \N$. Set the parameter $\lambda=q^n$.  Let $N_\lambda (a,b,c)$ denote the subspace of $M_\lambda(a,b,c)$ spanned by $\{m_i\}_{i=n+1}^\infty$. 
By construction $N_\lambda (a,b,c)$ is invariant under $A$. 
By (\ref{varphii}) the scalar $\varphi_{n+1}=0$. 
Combined with Theorem \ref{thm:Verma}(i) this implies that $N_\lambda (a,b,c)$ is invariant under $B$. 
By Theorem \ref{thm:Verma}(ii) the elements $\alpha,\beta,\gamma$ act on $N_\lambda (a,b,c)$ as scalar multiplication.
Thus $N_\lambda(a,b,c)$ is a $\triangle_q$-submodule of $M_\lambda(a,b,c)$ by Proposition \ref{prop:UAW_presentation}. 
Moreover $N_\lambda (a,b,c)$ is the $\triangle_q$-submodule of $M_\lambda(a,b,c)$ generated by $m_{n+1}$. 
Therefore
$$
V_n(a,b,c):=M_\lambda(a,b,c)/N_\lambda(a,b,c)
$$
is an $(n+1)$-dimensional $\triangle_q$-module that has the basis
$$
m_i+N_\lambda (a,b,c)
\qquad 
(0\leq i\leq n).
$$  
In this paper we will prove that all irreducible $\triangle_q$-modules with dimensions less than $\dbar$ are contained in the first family of finite-dimensional quotients of Verma $\triangle_q$-modules up to isomorphism:

\begin{thm}
\label{thm:<d}
Suppose that $V$ is an irreducible $\triangle_q$-module that has dimension less than $\dbar$. Then there exist an element $(a,b,c)\in {\F^\times}^3$ and an integer $n$ with $0\leq n\leq \dbar-2$
such that the $\triangle_q$-module $V_n(a,b,c)$ is isomorphic to $V$. 
\end{thm}

Let $\{\pm 1\}$ denote the multiplicative group consisting of the integers $1$ and $-1$.
There exists a unique left $\{\pm 1\}^3$-action on ${\F^{\times}}^3$ given by 
\begin{align*}
(-1,1,1)\cdot (a,b,c)&=(a^{-1},b,c),
\\
(1,-1,1)\cdot (a,b,c)&=(a,b^{-1},c),
\\
(1,1,-1)\cdot (a,b,c)&=(a,b,c^{-1})
\end{align*}
for all $(a,b,c)\in {\F^{\times}}^3$.
An irreducibility criterion for $V_n(a,b,c)$ can be expressed in terms of the $\{\pm 1\}^3$-action on ${\F^{\times}}^3$:

\begin{thm}
\label{thm:irrVn(abc)}
For any element $(a,b,c)\in {\F^\times}^3$ and any integer $n$ with $0\leq n\leq \dbar-2$ 
the following conditions are equivalent:
\begin{enumerate}
\item The $\triangle_q$-module $V_n(a,b,c)$ is irreducible.

\item $\bar{a}\bar{b}\bar{c}\not\in\{q^{n-2i+1}\,|\,i=1,2,\ldots,n\}$ for all $(\bar{a},\bar{b},\bar{c})\in \{\pm 1\}^3\cdot (a,b,c)$.
\end{enumerate}
\end{thm}
\begin{proof}
Similar to the proof of \cite[Theorem 4.4]{Huang:2015}.
\end{proof}

 Fix an integer $n$ with $0\leq n\leq \dbar-2$. 
Let $\I_n$ denote the set of the isomorphism classes of $(n+1)$-dimensional irreducible $\triangle_q$-modules. 
Let $\P_n$ denote the set of all elements $(a,b,c)\in {\F^\times}^3$ satisfying Theorem \ref{thm:irrVn(abc)}(ii). 
Clearly $\P_n$ is closed under the $\{\pm 1\}^3$-action on ${\F^{\times}}^3$. 
Similar to the case of $q$ not a root of unity, Theorems \ref{thm:<d} and \ref{thm:irrVn(abc)} result in a classification of all irreducible $\triangle_q$-modules that have dimensions less than $\dbar$:

\begin{thm}
For any integer $n$ with $0\leq n\leq \dbar-2$ there is a bijection $\{\pm 1\}^3 \backslash\P_n\to \I_n$ given by 
\begin{eqnarray*}
\{\pm 1\}^3\cdot (a,b,c)
&\mapsto &
\hbox{the isomorphism class of $V_n(a,b,c)$}
\end{eqnarray*}
for all $(a,b,c)\in \P_n$.
\end{thm}
\begin{proof}
Similar to the proof of \cite[Theorem 4.7]{Huang:2015}.
\end{proof}

We now take up the second family of finite-dimensional quotients of Verma $\triangle_q$-modules. 
Pick any quadruple $(a,b,c,\lambda)\in {\F^\times}^4$. 
Since $q^{2\dbar}=1$ the parameters (\ref{thetai})--(\ref{varphii}) satisfy  the cyclic properties:
\begin{align*}
\theta_i&=
\theta_{\dbar+i}
\qquad 
\hbox{for all $i\in \N$};
\\
\theta_i^*&=
\theta_{\dbar+i}^*
\qquad 
\hbox{for all $i\in \N$};
\\
\varphi_i &=
\varphi_{\dbar+i}
\qquad 
\hbox{for all $i\in \N$}.
\end{align*}  
Let $\delta\in \F$ denote an additional parameter. 
Define $O_{\lambda}^\delta(a,b,c)$ to be the subspace of $M_\lambda(a,b,c)$ spanned by 
$\{\delta m_i-m_{\dbar+i}\}_{i\in \N}$. 
Applying Theorem \ref{thm:Verma}(i) along with the cyclic properties, it is routine to verify that 
$O_\lambda^\delta(a,b,c)$ is invariant under $A$ and $B$.  
By Theorem \ref{thm:Verma}(ii) the elements $\alpha,\beta,\gamma$ act on $O_\lambda^\delta(a,b,c)$ as scalar multiplication.
Thus $O_\lambda^\delta(a,b,c)$ is a $\triangle_q$-submodule of $M_\lambda(a,b,c)$ by Proposition \ref{prop:UAW_presentation}. Moreover  $O_\lambda^\delta(a,b,c)$ is the $\triangle_q$-submodule of $M_\lambda(a,b,c)$ generated by 
$\delta m_0-m_{\dbar}$.
Therefore 
$$
W_\lambda^\delta (a,b,c):=M_\lambda(a,b,c)/O_\lambda^\delta(a,b,c)
$$
is a $\dbar$-dimensional $\triangle_q$-module that has the basis
\begin{gather*}
m_i+O_\lambda^\delta(a,b,c)
\qquad 
(0\leq i\leq \dbar-1).
\end{gather*}
In this paper we will prove that all $\dbar$-dimensional irreducible $\triangle_q$-modules with marginal weights are contained in the second family of finite-dimensional quotients of Verma $\triangle_q$-modules up to isomorphism:

\begin{thm}\label{thm:=d}
Suppose that $V$ is a $\dbar$-dimensional irreducible $\triangle_q$-module with marginal weights.
Then there exists an element $(a,b,c,\lambda,\delta)\in {\F^\times}^4\times \F$ such that the $\triangle_q$-module $W_{\lambda}^\delta(a,b,c)$ is isomorphic to $V$. 
\end{thm}

There is a unique left $\{\pm 1\}$-action on ${\F^{\times}}^4$ given by 
\begin{gather*}
(-1)\cdot (a,b,c,\lambda)
=
(-a,-b,-c,-\lambda)
%\label{z2s4-1}
\end{gather*}
for all $(a,b,c,\lambda)\in {\F^{\times}}^4$. 
Recall that the symmetric group $\S_4$ of degree four has a presentation given by the transpositions $(1\, 2), (2\, 3), (3\, 4)$ subject to the relations
\begin{gather*}
(1\, 2)^2=(2\, 3)^2=(3\, 4)^2=1,
\\
(1\, 2)(3\, 4)=(3\, 4)(1\, 2),
\\
(1\, 2)(2\, 3)(1\, 2)=(2\, 3)(1\, 2)(2\, 3),
\\
(2\, 3)(3\, 4)(2\, 3)=(3\, 4)(2\, 3)(3\, 4).
\end{gather*}
Applying the presentation it is straightforward to verify that there exists a unique right $\S_4$-action on $\{\pm 1\}\backslash{\F^{\times}}^4$ given by 
\begin{align*}
(\{\pm 1\}\cdot (a,b,c,\lambda))\cdot {(1\, 2)}
&=
\{\pm 1\}\cdot (a,b,c^{-1},\lambda),
%\label{z2s4-2}
\\
(\{\pm 1\}\cdot (a,b,c,\lambda))\cdot {(2\, 3)}
&=
\{\pm 1\}\cdot \textstyle (\frac{a}{\sqrt{a b c\lambda q}},
\frac{b}{\sqrt{a b c\lambda q}},
\frac{c}{\sqrt{a b c\lambda q}},
\frac{\lambda}{\sqrt{a b c\lambda q}}),
%\label{z2s4-3}
\\
(\{\pm 1\}\cdot (a,b,c,\lambda))\cdot {(3\, 4)}
&=
\{\pm 1\}\cdot (a^{-1},b,c,\lambda)
%\label{z2s4-4}
\end{align*}
for all $(a,b,c,\lambda)\in {\F^{\times}}^4$.

\begin{defn}
\label{defn:fz2s4&z2s4}
\begin{enumerate}
\item For any $(a,b,c,\lambda),(\bar{a},\bar{b},\bar{c},\bar{\lambda})\in  {\F^{\times}}^4$ we define 
$$
(a,b,c,\lambda)
\fz2s4\!
(\bar{a},\bar{b},\bar{c},\bar{\lambda})
$$ 
whenever $(\{\pm 1\}\cdot (a,b,c,\lambda))\cdot \S_4=(\{\pm 1\}\cdot (\bar{a},\bar{b},\bar{c},\bar{\lambda}))\cdot \S_4$ in $(\{\pm 1\}\backslash {\F^{\times}}^4)/ \S_4$.

\item 
For any $(a,b,c,\lambda,\delta),(\bar{a},\bar{b},\bar{c},\bar{\lambda},\bar{\delta})\in  {\F^{\times}}^4\times \F$ we define 
$$
(a,b,c,\lambda,\delta)
\z2s4\!
(\bar{a},\bar{b},\bar{c},\bar{\lambda},\bar{\delta})
$$ 
whenever $(a,b,c,\lambda)
\fz2s4\!
(\bar{a},\bar{b},\bar{c},\bar{\lambda})$ 
and 
\begin{gather}
\label{delta&bardelta}
\delta
+a^{\dbar}\lambda^{-\dbar}
+a^{-\dbar}\lambda^{\dbar}
=
\bar{\delta}
+\bar{a}^{\dbar}\bar{\lambda}^{-\dbar}
+\bar{a}^{-\dbar}\bar{\lambda}^{\dbar}.
\end{gather}
\end{enumerate}
Note that  
$\!\fz2s4$ is an equivalence relation on ${\F^{\times}}^4$ and $\!\z2s4$ is an equivalence relation on ${\F^{\times}}^4\times \F$.
\end{defn}

An irreducibility criterion for the $\triangle_q$-module $W_\lambda^\delta(a,b,c)$ can be expressed in terms of the equivalence relation $\!\z2s4$:

\begin{thm}
\label{thm2:irr}
For any $(a,b,c,\lambda,\delta)\in {\F^\times}^4\times \F$ the following conditions are equivalent:
\begin{enumerate}
\item The $\triangle_q$-module $W_\lambda^\delta(a,b,c)$ is irreducible.

\item 
$\bar{\delta}\not=0$ 
or 
$\bar{\lambda}^2\not\in\{q^{2i}\,|\,i=0,1,\ldots,\dbar-2\}$ 
for all 
$(\bar{a},\bar{b},\bar{c},\bar{\lambda},\bar{\delta})
\z2s4\!
(a,b,c,\lambda,\delta)$.
\end{enumerate}
\end{thm}

In order to classify the irreducible $\dbar$-dimensional $\triangle_q$-modules with marginal weights, the equivalence relation $\!\z2s4$ is extended as follows:

\begin{defn}
\label{defn:sim}
For any $
(a,b,c,\lambda,\delta), 
(\bar{a},\bar{b},\bar{c},\bar{\lambda},\bar{\delta})
\in {\F^{\times}}^4\times \F$ we define 
$
(a,b,c,\lambda,\delta)
\sim
(\bar{a},\bar{b},\bar{c},\bar{\lambda},\bar{\delta}) 
$ 
whenever any of the following conditions holds:
\begin{enumerate}
\item $(a,b,c,\lambda,\delta)\z2s4\!
(\bar{a},\bar{b},\bar{c},\bar{\lambda},\bar{\delta})$.

\item $(\bar{a},\bar{b},\bar{c},\bar{\lambda},\bar{\delta})=(a^{-1},b,c,\lambda^{-1}q^{-2},\delta)$ and $\lambda^2\in\{q^{2i}\,|\,i=0,1,\ldots,\dbar-2\}$.

\item 
$
(\bar{a},\bar{b},\bar{c},\bar{\lambda},\bar{\delta})
=
(a^{-1},b^{-1},c,\lambda^{-1}q^{-2},\delta)$ and the following conditions hold:
\begin{enumerate}
\item %$\lambda^2, a^{-1}b^{-1}c^{-1}\lambda q^{-1}, a^{-1}b^{-1}c\lambda q^{-1} \not\in \{q^{2i}\,|\,i=0,1,\ldots,\dbar-2\}$.
$b^2\lambda^{-2}\not\in \{q^{2(\dbar-i+1)}\,|\,i=0,1,\ldots,\dbar-2\}$.

\item 
$\delta (b^{\dbar}\lambda^{-\dbar}-b^{-\dbar}\lambda^{\dbar}) 
=a^{-\dbar}b^{-\dbar}
(\lambda^{2\dbar}-1)(a^{\dbar}b^{\dbar}c^{\dbar}\lambda^{-\dbar}q^{\dbar}-1)
(a^{\dbar}b^{\dbar}c^{-\dbar}\lambda^{-\dbar}q^{\dbar}-1)$.
\end{enumerate}
\end{enumerate}
\end{defn}

\begin{defn}
\label{defn:simeq}
Define $\simeq$ to be the equivalence relation on ${\F^\times}^4\times \F$ generated by $\sim$. 
\end{defn}

Let ${\bf IM}_{\dbar}$ denote the set of the isomorphism classes of all $\dbar$-dimensional irreducible $\triangle_q$-modules that have marginal weights.
Let ${\bf PM}_{\dbar}$ denote the set of all elements $(a,b,c,\lambda,\delta)\in {\F^\times}^4\times \F$ satisfying Theorem \ref{thm2:irr}(ii). 
It can be shown that ${\bf PM}_{\dbar}$ is closed under $\simeq$. 
Let ${\bf PM}_{\dbar}/\!\simeq$ denote the set of all equivalence classes of ${\bf PM}_{\dbar}$ under $\simeq$.

\begin{thm}
\label{thm:bijection}
There is a bijection from ${\bf PM}_{\dbar}/\!\simeq$ to ${\bf IM}_{\dbar}$ given by 
\begin{eqnarray*}
\hbox{the equivalence class of $(a,b,c,\lambda,\delta)$ under $\simeq$}
&\mapsto &
\hbox{the isomorphism class of $W_\lambda^\delta(a,b,c)$}
\end{eqnarray*}
for all $(a,b,c,\lambda,\delta)\in {\bf PM}_{\dbar}$.
\end{thm}

The outline of this paper is as follows: In Section \ref{s:marginal} we recall some properties of the marginal weights and closely related weight vectors called the marginal weight vectors. In Section \ref{s:feasible} we reinterpret the universal property of $M_\lambda(a,b,c)$ and relate the property to a functional relation called the feasible relation. In Section \ref{s:<d} we give a proof of Theorem \ref{thm:<d}. In Section \ref{s:polynomial} we establish a polynomial characterization for the feasible relation; consequently the feasible relation can be expressed in terms of the equivalence relation $\!\fz2s4$. In Section \ref{s:=d} we give a proof of Theorem \ref{thm:=d}. In Section \ref{s:fz2s4&z2s4} we characterize the equivalence relations $\!\fz2s4$ and $\!\z2s4$ in terms of the $\triangle_q$-modules $M_\lambda(a,b,c)$ and $W_\lambda^\delta(a,b,c)$. In Section \ref{s:irr} we give a proof of Theorem \ref{thm2:irr}. In Section \ref{s:sim} we relate the binary relation $\sim$ to the marginal weight vectors of $W_\lambda^\delta(a,b,c)$. In Section \ref{s:isoclass} we give a proof of Theorem \ref{thm:bijection}. In addition  the $\S_4$-action of $\{\pm 1\}\backslash {\F^{\times}}^4$ is fully displayed in Appendix \ref{s:S4}.

\section{The marginal weights and the marginal weight vectors}\label{s:marginal}

Recall the marginal weights of $\triangle_q$-modules from Definition \ref{defn:marginal}.

\begin{thm}
[Theorem 6.3, \cite{Huang:2021}]
\label{thm1:Huang2021}
The dimension of any finite-dimensional irreducible $\triangle_q$-module with marginal weights is less than or equal to $\dbar$.
\end{thm}

\begin{thm}
[Theorem 6.10, \cite{Huang:2021}]
\label{thm2:Huang2021}
The dimension of any finite-dimensional irreducible $\triangle_q$-module without marginal weights is equal to $\dbar$.
\end{thm}

\begin{lem}
\label{lem:<d&marginal}
Every nonzero $\triangle_q$-module with dimension less than $\dbar$ has marginal weights.
\end{lem}
\begin{proof}
Immediate from Theorems \ref{thm1:Huang2021} and \ref{thm2:Huang2021}. 
\end{proof}

\begin{defn}
\label{defn:marginalvector}
%Let $V$ denote a $\triangle_q$-module and let $\mu$ denote a weight of $V$.
Let $\mu$ denote a weight of a $\triangle_q$-module $V$. 
A weight vector $v$ of $V$ with weight $\mu$ is said to be {\it marginal} whenever $v$ is an eigenvector of 
$$
(B-\mu q^2-\mu^{-1} q^{-2}) A.
$$ 
\end{defn}

By Definition \ref{defn:marginal}, if a $\triangle_q$-module $V$ contains a marginal weight vector with weight $\mu$ then $\mu$ is a marginal weight of $V$.

\begin{lem}
[Lemma 6.1, \cite{Huang:2021}]
\label{lem1:marginalvector}
Assume that $V$ is a finite-dimensional irreducible $\triangle_q$-module. For any weight $\mu$ of $V$ the following conditions are equivalent:
\begin{enumerate}
\item $\mu$ is a marginal weight of $V$.

\item There exists a marginal weight vector of $V$ with weight $\mu$.
\end{enumerate}
\end{lem}

\begin{lem}
\label{lem2:marginalvector}
Assume a finite-dimension irreducible $\triangle_q$-module $V$ contains a marginal weight vector $v$ with weight $\mu$. 
For all $i\in \N$ the following statements hold:
\begin{enumerate}
\item $
(B-\mu q^{2i}-\mu^{-1} q^{-2i})A^i v
$ is a linear combination of $v, Av,\ldots, A^{i-1} v$. 

\item $
\prod\limits_{h=0}^i (B-\mu q^{2h}-\mu^{-1} q^{-2h})
$ vanishes at $v, Av,\ldots, A^i v$.
\end{enumerate}
\end{lem}
\begin{proof}
(i): Immediate from \cite[Lemma 6.2]{Huang:2021}.

(ii): Applying Lemma \ref{lem2:marginalvector}(i) a routine induction on $i$ yields (ii).
\end{proof}

\begin{lem}
\label{lem3:marginalvector}
If a finite-dimensional irreducible $\triangle_q$-module $V$ contains a marginal weight vector $v$, then $V$ is spanned by $\{A^i v\}_{i\in \N}$.
\end{lem}
\begin{proof}
Let $W$ denote the subspace of $V$ spanned by $\{A^i v\}_{i\in \N}$. Then $W$ is $A$-invariant. By Lemma \ref{lem2:marginalvector}(i), $W$ is $B$-invariant. By Schur's lemma the central elements $\alpha,\beta,\gamma$ act on $V$, as well as $W$, as scalar multiplication. Hence $W$ is a $\triangle_q$-submodule of $V$ by Proposition \ref{prop:UAW_presentation}. Since $v\in W$ and $v\not=0$ it follows that $W=V$. The lemma follows.
\end{proof}

\section{The universal property of $M_\lambda(a,b,c)$ and the feasible relation}\label{s:feasible}

For the sake of brevity, the following notational agreements will be used throughout the rest of this paper. 
While using the quadruple $(a,b,c,\lambda)$ to represent an element of ${\F^\times}^4$, the notation $\{m_i\}_{i\in \N}$ denotes the basis for $M_\lambda(a,b,c)$ mentioned in Theorem \ref{thm:Verma}(i) and $\{\theta_i\}_{i\in \N}$, $\{\theta_i^*\}_{i\in \N}$, $\{\varphi_i\}_{i\in \N}$ always stand for the accompanying parameters (\ref{thetai})--(\ref{varphii}).

We begin with this section with a simplified description for \cite[Proposition 3.2]{Huang:2015}, which is called the universal property of the Verma $\triangle_q$-module $M_\lambda(a,b,c)$.

\begin{prop}
\label{prop2:Verma_universal}
Let $(a,b,c,\lambda)\in {\F^\times}^4$ be given. 
For any $\triangle_q$-module $V$ and $v\in V$
the following conditions are equivalent:
\begin{enumerate}
\item There exists a $\triangle_q$-module homomorphism 
$
M_\lambda(a,b,c)\to V
$ 
that sends $m_0$ to $v$.

\item The following equations hold on $V$:  
\begin{gather}
B v=\theta_0^* v,
\label{universal'-1}
\\
(B-\theta_1^*)A v=(\theta_0(\theta_0^*-\theta_1^*)+\varphi_1)  v,
\label{universal'-2}
\\
\beta v=\omega^* v,
\label{universal'-3}
\\
\gamma v=\omega^\e v,
\label{universal'-4}
\end{gather}
where $\omega^*$ and $\omega^\e$ are the scalars {\rm (\ref{omega*})} and {\rm (\ref{omegae})}, respectively. 
\end{enumerate}
\end{prop}
\begin{proof}
By Theorem \ref{thm:Verma} the condition (i) implies 
(\ref{universal'-1}), (\ref{universal'-3}), (\ref{universal'-4}) and the following equations:
\begin{gather}
(B-\theta_1^*)(A-\theta_0) v=\varphi_1  v,
\label{universal-2}
\\
\alpha v=\omega v, 
\label{universal-3}
\end{gather}
where $\omega$ is the scalar (\ref{omega}). 
By \cite[Proposition 3.2]{Huang:2015} the equations (\ref{universal'-1}) and (\ref{universal'-3})--(\ref{universal-3}) imply the condition (i). Observe that (\ref{universal'-2}) is identical to (\ref{universal-2}) when (\ref{universal'-1}) holds. 
Applying $v$ to either side of (\ref{e:a->ABc}),  
we evaluate the resulting equation by using (\ref{universal'-1}), (\ref{universal'-2}) and (\ref{universal'-4}); thereby gaining the equation (\ref{universal-3}). 
Therefore (\ref{universal'-1})--(\ref{universal'-4}) hold if and only if (\ref{universal'-1}) and (\ref{universal'-3})--(\ref{universal-3}) hold.
The proposition follows. 
\end{proof}

Since the $\triangle_q$-module $M_\lambda(a,b,c)$ is generated by $m_0$, if Proposition \ref{prop2:Verma_universal}(i) holds then the mentioned map is unique. The coefficient of $v$ in the right-hand side of (\ref{universal'-2}) is equal to $q-q^{-1}$ times 
$$
(c+c^{-1})(\lambda-\lambda^{-1})-(a+a^{-1})(b q-b^{-1} q^{-1}).
$$
Inspired by Proposition \ref{prop2:Verma_universal} we study the following functional relation:

\begin{defn}
\label{defn:feasible}
For any $(a,b,c,\lambda)\in \F^{\times 4}$ and any $(\mu,\varphi,\omega^*,\omega^\e)\in \F^\times \times \F^3$, we say that $(a,b,c,\lambda)$ is {\it feasible} for $(\mu,\varphi,\omega^*,\omega^\e)$ whenever the following equations hold:
\begin{enumerate}
\item $\mu=b\lambda^{-1}$.

\item $\varphi=
(c+c^{-1})(\lambda-\lambda^{-1})-(a+a^{-1})(b q-b^{-1} q^{-1})$.

\item $\omega^{*}=
(c+c^{-1})(a+a^{-1})+(b+b^{-1})(\lambda q+\lambda^{-1} q^{-1})$.

\item $\omega^{\e}=(a+a^{-1})(b+b^{-1})+(c+c^{-1})(\lambda q+\lambda^{-1} q^{-1})$.
\end{enumerate}
\end{defn}

\begin{thm}
\label{thm:feasible&universal}
Let $(a,b,c,\lambda)\in \F^{\times 4}$ and $(\mu,\varphi,\omega^*,\omega^\e)\in \F^\times \times \F^3$ be given. 
Suppose that a $\triangle_q$-module $V$ contains a nonzero vector $v$ satisfying the following equations:
\begin{gather}
Bv=(\mu+\mu^{-1}) v,
\label{f&u-1}
\\
(B-\mu q^2-\mu^{-1}q^{-2}) A v=(q-q^{-1})\varphi\cdot v,
\label{f&u-2}
\\
\beta v=\omega^{*}v,
\label{f&u-3}
\\
\gamma v=\omega^{\e} v.
\label{f&u-4}
\end{gather}
Then $(a,b,c,\lambda)$ is feasible for $(\mu,\varphi,\omega^*,\omega^\e)$ if and only if the following conditions hold:
\begin{enumerate}
\item $\mu=b\lambda^{-1}$.

\item There exists a $\triangle_q$-module homomorphism 
$
M_\lambda(a,b,c)\to V
$ 
that sends $m_0$ to $v$. 
\end{enumerate}
\end{thm}
\begin{proof}
The condition (i) is exactly
Definition \ref{defn:feasible}(i). 

($\Rightarrow$): 
The condition (ii) is immediate from Proposition \ref{prop2:Verma_universal} and Definition \ref{defn:feasible}.

($\Leftarrow$): 
Since (ii) holds it follows from Proposition \ref{prop2:Verma_universal} that the equations (\ref{universal'-1})--(\ref{universal'-4}) hold. By (i) the scalar $\theta_1^*$ in  (\ref{universal'-3}) is equal to $\mu q^2+\mu^{-1}q^{-2}$. 
Since $v\not=0$ and comparing (\ref{universal'-2}) with (\ref{f&u-2}), this yields Definition \ref{defn:feasible}(ii). 
For similar reasons 
Definition \ref{defn:feasible}(iii), (iv)
follow. Therefore $(a,b,c,\lambda)$ is feasible for $(\mu,\varphi,\omega^*,\omega^\e)$.
\end{proof}

We end this section with the following lemmas related to Theorem \ref{thm:feasible&universal}:

\begin{lem}
\label{lem:irr&feasible}
Suppose that $V$ is a finite-dimensional irreducible $\triangle_q$-module with a marginal weight $\mu$. Then there are a nonzero vector $v\in V$ and three scalars $\varphi,\omega^*,\omega^\e\in \F$ satisfying the equations {\rm (\ref{f&u-1})--(\ref{f&u-4})}.
\end{lem}
\begin{proof}
By Lemma \ref{lem1:marginalvector} there exists a marginal weight vector $v$ of $V$ with weight $\mu$. Hence (\ref{f&u-1}) follows. 
By Definition \ref{defn:marginalvector} and since $q^2\not=1$ there is a scalar $\varphi\in \F$ such that (\ref{f&u-2}) holds. 
By Schur's lemma there are $\omega^*,\omega^\e\in \F$ such that (\ref{f&u-3}) and (\ref{f&u-4}) hold. 
\end{proof}

\begin{lem}
\label{lem:Verma_marginalweight}
For any $\triangle_q$-submodule $O$ of $M_\lambda(a,b,c)$ with $m_0\not\in O$ the following statements are true:
\begin{enumerate}
\item $m_0+O$ is a marginal weight vector of $M_\lambda(a,b,c)/O$ with weight $b\lambda^{-1}$.

\item $m_0+O$ is a marginal weight vector of $M_\lambda(a,b,c)/O$ with weight $b^{-1}\lambda$ if and only if 
$b\lambda^{-1}=b^{-1}\lambda$, or $m_0+O$ and $m_1+O$ are linearly dependent.
\end{enumerate}
\end{lem}
\begin{proof}
(i): By Theorem \ref{thm:Verma}(i) it is routine to verify the statement (i).

(ii): By Theorem \ref{thm:Verma}(i) a straightforward calculation shows that 
\begin{align*}
(B-b^{-1}\lambda q^2-b\lambda^{-1}q^{-2}) 
Am_0=
%(q-q^{-1})
\varphi\cdot 
m_0
+
(q^2-q^{-2})(b\lambda^{-1}-b^{-1}\lambda) m_1
\end{align*}
for some scalar $\varphi\in \F$.
%where 
%$
%\varphi=(q+q^{-1}) (ab\lambda^{-2}-a^{-1}b^{-1}\lambda^2)
%+(c+c^{-1})(\lambda-\lambda^{-1})
%-(b+b^{-1})(a q-a^{-1} q^{-1})$. 
The statement (ii) follows from the above equation.
\end{proof}

\section{Proof for Theorem \ref{thm:<d}}
\label{s:<d}

\begin{lem}
\label{lem:vee}
There exists a unique algebra automorphism $\vee$ of $\triangle_q$ that sends 
\begin{eqnarray*}
(A,B,\alpha,\beta,\gamma) &\mapsto & 
(B,A,\beta,\alpha,\gamma).
\end{eqnarray*}
\end{lem}
\begin{proof}
It is routine to verify the lemma by using Proposition \ref{prop:UAW_presentation} .
\end{proof}

For any $\triangle_q$-module $V$ the notation 
$$
V^\vee
$$ 
stands for the $\triangle_q$-module obtained  by twisting the $\triangle_q$-module $V$ via $\vee$.

\begin{lem}
\label{lem:<d_Deltai}
Let $n$ denote an integer with $1\leq n\leq \dbar-2$. Suppose that $\{h_i\}_{i=0}^{n+1}$ is a sequence in $\F$ satisfying the following conditions:
\begin{enumerate}
\item $h_{i+2}
-h_{i-1}
=
(q^2+1+q^{-2})(h_{i+1}-h_i)
$ 
for all $i=1,2,\ldots,n-1$.

\item 
There are three integers $j,k,\ell$ with $0\leq j<k<\ell\leq n+1$ such that $h_j=h_k=h_\ell=0$.
\end{enumerate}
Then $
h_i=0$ for all $i=0,1,\ldots,n+1$.
\end{lem}
\begin{proof}
If $n=1$ then the lemma is immediate from (ii).
Now suppose that $n\geq 2$. Since $q^4\not=1$ the roots $1$, $q^2$, $q^{-2}$ of the characteristic equation for the linear homogeneous recurrence (i) are mutually distinct. Thus there are $c_0,c_1,c_2\in \F$ such that 
$$
h_i=c_0 + q^{2i} c_1 + q^{-2i} c_2 
\qquad 
(0\leq i\leq n+1).
$$
By (ii) the following linear equations hold:
\begin{align*}
\left\{
\begin{array}{ll}
c_0 +  q^{2j}c_1 + q^{-2j} c_2 =0,
\\
c_0 +  q^{2k}c_1 + q^{-2k} c_2 =0,
\\
c_0 +  q^{2\ell}c_1 + q^{-2\ell} c_2  =0.
\end{array}
\right.
\end{align*}
The determinant of the coefficient matrix for the above linear equations is equal to 
$$
(q^{j-k}-q^{k-j})(q^{k-\ell}-q^{\ell-k})(q^{\ell-j}-q^{j-\ell}).
$$
Since each of $k-j$, $\ell-k$, $\ell-j$ is a positive integer less than $\dbar$, none of $q^{2(k-j)}$, $q^{2(\ell-k)}$, $q^{2(\ell-j)}$ is equal to one. Therefore the determinant is nonzero. Since the coefficient matrix is invertible each of $c_0, c_1, c_2$ is zero. The lemma follows.
\end{proof}

%\medskip

\begin{proof}[Proof of Theorem \ref{thm:<d}]
%\noindent{\it Proof of Theorem \ref{thm:<d}.} 
By Lemma \ref{lem:<d&marginal} 
there exists a marginal weight $\mu$ of the $\triangle_q$-module $V$. By Lemma \ref{lem:irr&feasible} there are a nonzero vector $v$ of $V$ and three scalars $\varphi,\omega^*,\omega^\e\in \F$ satisfying the equations (\ref{f&u-1})--(\ref{f&u-4}).

By Lemma \ref{lem:<d&marginal} there exists a marginal weight $\kappa$ of the $\triangle_q$-module $V^\vee$. Let $n$ denote the dimension of $V$ minus one. 
Set
\begin{gather}\label{<d:ablambda}
(a,b,\lambda)=(\kappa q^n,\mu q^n,q^n)
\end{gather}
and $c$ to be the scalar in $\F^\times$ satisfying
\begin{gather}\label{<d:c}
c+c^{-1}=\left\{
\begin{array}{ll}
\displaystyle 
\frac{\omega^\e-(a+a^{-1})(b+b^{-1})}{\lambda q+\lambda^{-1}q^{-1}}
\qquad 
&\hbox{if $n=0$},
\\
\displaystyle 
\frac{\varphi+(a+a^{-1})(bq-b^{-1}q^{-1})}{\lambda-\lambda^{-1}}
\qquad 
&\hbox{if $n\geq 1$}.
\end{array}
\right.
\end{gather}
Since $0\leq n< \dbar$ it follows that $\lambda^2=1$ if and only if $n=0$. Since $q^4\not=1$ it follows that $\lambda^2 \not=-q^{-2}$ when $n=0$. Hence the denominators in the right-hand side of (\ref{<d:c}) are nonzero. Since $\F$ is algebraically closed the existence of $c$ follows.

We are now going to show that $(a,b,c,\lambda)$ is feasible for $(\mu,\varphi,\omega^*,\omega^\e)$. 
Apparently Definition \ref{defn:feasible}(i) is immediate from (\ref{<d:ablambda}). 
Due to (\ref{<d:c}) we divide the argument into the two cases: $n=0$ and $n\geq 1$.

$(n=0)$: 
In this case $v$ is a basis for $V$ and $(a,b,\lambda)=(\kappa,\mu,1)$. 
Definition \ref{defn:feasible}(iv) is immediate from (\ref{<d:c}). 
Since $a$ is a weight of $V^\vee$ it follows from Lemma \ref{lem:vee} that 
\begin{align}
\label{n=0:A}
Av=(a+a^{-1}) v
\end{align} 
on $V$.
Evaluating the left-hand side of (\ref{f&u-2}) by using (\ref{f&u-1}) and (\ref{n=0:A}) yields that 
$$
\varphi=-(a+a^{-1})(b q-b^{-1}q^{-1}).
$$
Definition \ref{defn:feasible}(ii) follows. Applying $v$ to either side of (\ref{e:b->ABc}) we evaluate the resulting equation by using (\ref{f&u-1}) and (\ref{n=0:A}). It follows that $\omega^*$ is equal to 
$$
(a+a^{-1})
\cdot \frac{\omega^\e-(a+a^{-1})(b+b^{-1})}{q+q^{-1}}
+
(q+q^{-1})(b+b^{-1}).
$$
Definition \ref{defn:feasible}(iii) follows by substituting Definition \ref{defn:feasible}(iv) into the above scalar. 
Therefore $(a,b,c,\lambda)$ is feasible for $(\mu,\varphi,\omega^*,\omega^\e)$ when $n=0$.

$(n\geq 1)$: 
Definition \ref{defn:feasible}(ii) is immediate from (\ref{<d:c}).
%Let $\{\theta_i\}_{i\in \N},\{\theta_i^*\}_{i\in \N},\{\varphi_i\}_{i\in \N}$ denote the parameters (\ref{thetai})--(\ref{varphii}) associated with $a,b,c,\lambda$ given by (\ref{<d:ablambda}) and (\ref{<d:c}).
Set
\begin{gather}\label{<d:vi}
v_0=v,
\qquad 
v_i=(A-\theta_{i-1})  v_{i-1}
\qquad 
(1\leq i\leq n).
\end{gather}
By Lemma \ref{lem3:marginalvector} the vectors $\{v_i\}_{i=0}^n$ are a basis for $V$.  
By Lemma \ref{lem1:marginalvector} there exists a marginal weight vector $w$ of the $\triangle_q$-module  $V^\vee$ with weight $\kappa$. 
It follows from Lemma \ref{lem2:marginalvector}(ii) that 
$
\prod_{i=0}^n (A-\theta_i)
$
vanishes at $\{B^i w\}_{i=0}^n$ on the $\triangle_q$-module $V$. By Lemma \ref{lem3:marginalvector}
the vectors $\{B^i w\}_{i=0}^n$ are also a basis for $V$. 
Hence $
\prod_{i=0}^n (A-\theta_i)
$
vanishes on $V$. 
Combined with (\ref{<d:vi}) this implies that 
\begin{gather}\label{<d:vn}
(A-\theta_n)v_n=0.
\end{gather}
Let $\varphi^{(i)}_j\in \F$ for all integers $i,j$ with $0\leq i,j\leq n$ such that 
\begin{gather}\label{<d:B}
(B-\theta_i^*)v_i=\sum_{j=0}^n \varphi^{(i)}_j v_j.
\end{gather}
By Lemma \ref{lem2:marginalvector}(i) the coefficients 
\begin{gather}\label{<d:varphij=0}
\varphi^{(i)}_j=0 
\qquad 
(0\leq i\leq j\leq n).
\end{gather}

Applying (\ref{f&u-1}) and (\ref{f&u-2}) yields that 
$$
\varphi^{(1)}_0=(q-q^{-1})\varphi-\theta_0(\theta_0^*-\theta_1^*).
$$  
Evaluating the right-hand side of the above equation by using 
Definition \ref{defn:feasible}(ii)
yields that 
\begin{gather}\label{<d:varphi1}
\varphi^{(1)}_0=\varphi_1.
\end{gather}
For any integer $i$ with $1\leq i\leq n$ we apply $v_{i-1}$ to either side of (\ref{e:b->ABc}) and evaluate the coefficient of $v_i$ in the resulting equation by using (\ref{f&u-1}) and (\ref{<d:vi})--(\ref{<d:varphij=0}). It follows that 
\begin{equation}
\label{<d:omegae}
%\varphi_{i+1}'
%-(q^2+q^{-2})\varphi_i'
%+\varphi_{i-1}'
\varphi^{(i+1)}_i
-(q^2+q^{-2}) \varphi^{(i)}_{i-1}
+\varphi^{(i-1)}_{i-2}
+c_i
+(q-q^{-1})^2\omega^\e
=0
\qquad (1\leq i\leq n).
\end{equation}
Here 
$\varphi^{(0)}_{-1}=\varphi^{(n+1)}_n=0$ and 
$$
c_i=(\theta_i+\theta_{i-1})(\theta_i^*+\theta_{i-1}^*)-(q^2+q^{-2})(\theta_i\theta_i^*+\theta_{i-1}\theta_{i-1}^*)
\qquad 
(1\leq i\leq n).
$$
%for all integers $i$ with $1\leq i\leq n$.
Applying (\ref{<d:omegae}) yields that 
\begin{gather*}%\label{<d:varphi'_rec}
\varphi^{(i+2)}_{i+1}
-(q^2+1+q^{-2})
(\varphi^{(i+1)}_i-\varphi^{(i)}_{i-1})
-\varphi^{(i-1)}_{i-2}
=c_i-c_{i+1}
\qquad 
(1\leq i\leq n-1).
\end{gather*}
On the other hand a direct calculation shows that $\{\varphi_i\}_{i=0}^{n+1}$ also satisfy the above recurrence. 
Hence the scalars 
$$
h_i=\varphi^{(i)}_{i-1}-\varphi_i
\qquad 
(0\leq i\leq n+1)
$$  
satisfy Lemma \ref{lem:<d_Deltai}(i). 
By (\ref{<d:varphi1}) and since $\varphi_0=\varphi_{n+1}=\varphi^{(0)}_{-1}=\varphi^{(n+1)}_n=0$ it follows that $h_0=h_1=h_{n+1}=0$. 
%Hence $\{h_i\}_{i=0}^{n+1}$ satisfy Lemma \ref{lem:<d_Deltai}(ii). 
Applying Lemma \ref{lem:<d_Deltai} yields that 
\begin{gather*} %\label{<d:varphi}
\varphi^{(i)}_{i-1}=\varphi_i
\qquad 
(0\leq i\leq n+1).
\end{gather*}
It is now routine to verify 
Definition \ref{defn:feasible}(iv) 
by substituting the above equations into (\ref{<d:omegae}).

For any integer $i$ with $0\leq i\leq n$ 
we apply $v_i$ to either side of (\ref{e:b->ABc}) and evaluate the coefficient of $v_i$ in the resulting equation. It follows that 
\begin{align}\label{<d:varphi''_rec}
(q-q^{-1})(q^2-q^{-2})\omega^*=
%\varphi_{i+2}''
\varphi^{(i+2)}_i
-(q^2+q^{-2})
%\varphi_{i+1}''
\varphi^{(i+1)}_{i-1}
+
%\varphi_{i}''
\varphi^{(i)}_{i-2}
+ 
c_i
\qquad 
(0\leq i\leq n).
\end{align}
Here 
$\varphi^{(0)}_{-2}=\varphi^{(1)}_{-1}=\varphi^{(n+1)}_{n-1}=\varphi^{(n+2)}_n=0$ and 
\begin{align*}
\begin{split}
c_i
&=
%(\theta_i+\theta_{i-1}) \varphi_i
%-(q^2+q^{-2})\theta_i(\varphi_i+\varphi_{i+1})
%+(\theta_i+\theta_{i+1})\varphi_{i+1}
(\theta_i-\theta_{i+1})\varphi_i
+(\theta_i-\theta_{i-1})\varphi_{i+1}
-(q-q^{-1})^2\theta_i^2\theta_i^*
\\
&
\qquad
+\,(q^2-q^{-2})^2\theta_i^*
+(q-q^{-1})^2\theta_i\omega^\e
\end{split}
\qquad 
(0\leq i\leq n).
\end{align*}
A direct calculation shows that $c_i$ is equal to $(q-q^{-1})(q^2-q^{-2})$ times the right-hand side of 
Definition \ref{defn:feasible}(iii) for each $i=0,1,\ldots,n$. 
Combined with (\ref{<d:varphi''_rec}) this implies that 
$
\varphi^{(i+2)}_i
-(q^2+q^{-2})\varphi^{(i+1)}_{i-1}
+\varphi^{(i)}_{i-2}
$ 
for all $i=0,1,\ldots,n$ 
are identical. It follows that 
$$
\varphi^{(i+2)}_{i}-(q^2+1+q^{-2})(\varphi^{(i+1)}_{i-1}-\varphi^{(i)}_{i-2})-\varphi^{(i-1)}_{i-3}=0
\qquad 
(1\leq i\leq n).
$$
Applying Lemma \ref{lem:<d_Deltai} yields that 
$$
\varphi^{(i)}_{i-2}=0 
\qquad (0\leq i\leq n+2). 
$$
Definition \ref{defn:feasible}(iii) follows by substituting the above equations into (\ref{<d:varphi''_rec}). Therefore $(a,b,c,\lambda)$ is feasible for $(\mu,\varphi,\omega^*,\omega^\e)$ when $n\geq 1$.

Since $(a,b,c,\lambda)$ is feasible for $(\mu,\varphi,\omega^*,\omega^\e)$ it follows from Theorem \ref{thm:feasible&universal} that there exists a $\triangle_q$-module homomorphism 
\begin{gather}\label{<d:M->V}
M_\lambda(a,b,c)\to V
\end{gather}
that sends $m_0$ to $v$. By (\ref{n=0:A}) and (\ref{<d:vn}) the vector $m_{n+1}$ lies in the kernel of (\ref{<d:M->V}). Recall from Section \ref{s:intro} that the $\triangle_q$-submodule $N_\lambda(a,b,c)$ of $M_\lambda(a,b,c)$ is generated by $m_{n+1}$. Hence (\ref{<d:M->V}) induces the $\triangle_q$-module homomorphism 
\begin{gather}
\label{<d:Vn(abc)->V}
V_n(a,b,c)\to V
\end{gather}
that sends $m_0+N_\lambda(a,b,c)$ to $v$. Since $V_n(a,b,c)$ and $V$ have the same dimension $n+1$
and by the irreducibility of the $\triangle_q$-module $V$, the map (\ref{<d:Vn(abc)->V}) is a $\triangle_q$-module isomorphism. The result follows. 
%\hfill $\square$
\end{proof}

\section{A polynomial characterization for the feasible relation}\label{s:polynomial}

Recall the feasible relation from Definition \ref{defn:feasible}.
We derive the following characterization for the feasible relation:

\begin{thm}
\label{thm:universal_class}
%Assume that $(\kappa,c,\lambda)\in {\F}^{\times 3}$ and 
Assume that $(\mu,\varphi,\omega^*,\omega^\e)\in \F^\times \times \F^3$. For any scalars $\kappa,\lambda,c\in \F^\times$ the quadruple 
$
(\kappa\lambda,\mu\lambda,c,\lambda)
$ 
is feasible for $(\mu,\varphi,\omega^*,\omega^\e)$
if and only if the following conditions hold:
\begin{enumerate}
\item $\kappa$ is a root of the polynomial
$$
\frac{x^4}{\mu q}
-\frac{\omega^\e+q^{-1}\varphi}{q+q^{-1}}x^3
+(\omega^*-\mu q^{-1}-\mu^{-1} q)x^2
-\frac{\omega^\e-q\varphi}{q+q^{-1}}x
+\mu q.
$$

\item $\lambda$ is a root of the polynomial 
$$
\kappa\mu q x^6
+
\left(
\kappa^{-1}\mu q
-\frac{\omega^\e-q\varphi}{q+q^{-1}}
\right)x^4
+
\left(
\frac{\omega^\e+q^{-1}\varphi}{q+q^{-1}}
-
\kappa\mu^{-1} q^{-1}
\right)x^2
-\frac{1}{\kappa\mu q}.
$$

\item $c$ is a root of $x^2-rx+1$ where
\begin{align*}
r=
\left\{
\begin{array}{ll}
\displaystyle 
\frac{\varphi+(\kappa\lambda +\kappa^{-1}\lambda^{-1})(\mu\lambda q-\mu^{-1}\lambda^{-1}q^{-1})}{\lambda-\lambda^{-1}}
\qquad 
&\hbox{if $\lambda^2\not=1$},
\\
\displaystyle 
\frac{\omega^\e-(\kappa \lambda+\kappa^{-1}\lambda^{-1})(\mu\lambda+\mu^{-1}\lambda^{-1})}{\lambda q+\lambda^{-1} q^{-1}}
\qquad 
&\hbox{if $\lambda^2=1$}.
\end{array}
\right.
\end{align*}
\end{enumerate}
\end{thm}
\begin{proof}
Let 
\begin{gather}
\label{poly:ab}
(a,b)=(\kappa\lambda,\mu\lambda).
\end{gather}
Clearly  
Definition \ref{defn:feasible}(i) 
holds. Then $(a,b,c,\lambda)$ is feasible for $(\mu,\varphi,\omega^*,\omega^\e)$ if and only if 
Definition \ref{defn:feasible}(ii)--(iv)
hold.

($\Rightarrow$): 
Suppose that Definition \ref{defn:feasible}(ii)--(iv) hold and we show (i)--(iii). The condition (iii) is immediate from 
Definition \ref{defn:feasible}(ii) if $\lambda^2\not=1$; 
the condition (iii) is immediate from 
Definition \ref{defn:feasible}(iv) if $\lambda^2=1$.

Using Definition \ref{defn:feasible}(ii), (iv) it is routine to verify that  
\begin{align}
\frac{\omega^\e-q\varphi}{q+q^{-1}}
&=\lambda^{-1}(c+c^{-1})
+bq(a+a^{-1}),
\label{alg:omegae-qphi}
\\
\frac{\omega^\e+q^{-1}\varphi}{q+q^{-1}}
&=\lambda(c+c^{-1})
+b^{-1}q^{-1}(a+a^{-1}).
\label{alg:omegae+q-1phi}
\end{align}
We multiply (\ref{alg:omegae-qphi}) and (\ref{alg:omegae+q-1phi}) by $\lambda$ and $\lambda^{-1}$ respectively. The difference of the resulting equations gives 
$$
\lambda\frac{\omega^\e-q\varphi}{q+q^{-1}}-\lambda^{-1}\frac{\omega^\e+q^{-1}\varphi}{q+q^{-1}}
=(a+a^{-1})(b\lambda q-b^{-1} \lambda^{-1} q^{-1}).
$$
The condition (ii) follows by substituting (\ref{poly:ab}) into the above equation.

We multiply (\ref{alg:omegae-qphi}) and (\ref{alg:omegae+q-1phi}) by $a^{-1}\lambda$ and $a\lambda^{-1}$ respectively. The sum of the resulting equations gives
\begin{gather}
\label{alg:omegae-qphi+omegae+q-1phi}
a^{-1}\lambda \frac{\omega^\e-q\varphi}{q+q^{-1}}+a\lambda^{-1} \frac{\omega^\e+q^{-1}\varphi}{q+q^{-1}}
=
(a+a^{-1})(a^{-1}b\lambda q+ab^{-1}\lambda^{-1}q^{-1}+c+c^{-1}).
\end{gather}
%Using Definition \ref{defn:feasible}(iii) yields that 
%\begin{gather}
%\label{alg:omega*-muq-1-mu-1q}
%\omega^*-b\lambda^{-1}q^{-1}-b^{-1}\lambda q
%=(a+a^{-1})(c+c^{-1})+b\lambda q+b^{-1}\lambda^{-1} q^{-1}.
%\end{gather}
Subtracting Definition \ref{defn:feasible}(iii) from (\ref{alg:omegae-qphi+omegae+q-1phi}) yields that 
$$
a^{-1}\lambda 
\frac{\omega^\e-q\varphi}{q+q^{-1}}
+a\lambda^{-1} 
\frac{\omega^\e+q^{-1}\varphi}{q+q^{-1}}
-\omega^*
%+
%b\lambda^{-1}q^{-1}+b^{-1}\lambda q
=
%a^2b^{-1}\lambda^{-1}q^{-1}+a^{-2}b\lambda q.
(ab^{-1}-a^{-1}b)(a\lambda^{-1}q^{-1}-a^{-1}\lambda q).
$$
The condition (i) follows by substituting (\ref{poly:ab}) into the above equation. The ``only if'' part follows.

($\Leftarrow$): 
Suppose that (i)--(iii) hold and we show  
Definition \ref{defn:feasible}(ii)--(iv). 
Definition \ref{defn:feasible}(ii) is immediate from (iii) if $\lambda^2\not=1$; Definition \ref{defn:feasible}(ii) is obtained  from (ii) if $\lambda^2=1$.

Definition \ref{defn:feasible}(iv)
is immediate from (iii) if $\lambda^2=1$.
Suppose that $\lambda^2\not=1$. Applying (ii) yields that 
$$
\frac{\lambda-\lambda^{-1}}{q+q^{-1}}\omega^\e=(\kappa\lambda+\kappa^{-1}\lambda^{-1})(\mu\lambda^2 q-\mu^{-1}\lambda^{-2} q^{-1})
+\frac{\lambda q+\lambda^{-1} q^{-1}}{q+q^{-1}}\varphi.
$$ 
To get 
Definition \ref{defn:feasible}(iv), 
we evaluate the right-hand side of the above equation by using 
Definition \ref{defn:feasible}(ii) and replacing $\kappa$ and $\mu$ with $a\lambda^{-1}$ and $b\lambda^{-1}$.

Applying (i) yields that
$$
\omega^*=(\kappa^{-1}q-\kappa q^{-1})
\left(
\kappa\mu^{-1}-\kappa^{-1}\mu-\frac{\varphi}{q+q^{-1}}
\right)
+\frac{\kappa+\kappa^{-1}}{q+q^{-1}}\omega^\e.
$$
To get Definition \ref{defn:feasible}(iii),
we evaluate the right-hand side of the above equation by using 
Definition \ref{defn:feasible}(ii), (iv)
and replacing $\kappa$ and $\mu$ with $a\lambda^{-1}$ and $b\lambda^{-1}$. The ``if'' part follows.
\end{proof}

By Theorem \ref{thm:universal_class} there are at most $48$ elements of ${\F^\times}^4$ which are feasible for a given element of $\F^\times \times \F^{3}$. By Definition \ref{defn:fz2s4&z2s4}(i) each equivalence class under $\!\fz2s4$ consists of at most $48$ elements.

\begin{thm}
\label{thm:feasible}
For any $(\mu,\varphi,\omega^*,\omega^\e)\in \F^\times \times \F^{3}$ the following statements are true:
\begin{enumerate}
\item There exists an element $(a,b,c,\lambda)\in {\F^\times}^4$ which is feasible for $(\mu,\varphi,\omega^*,\omega^\e)$.

\item If $(a,b,c,\lambda)$ is feasible for $(\mu,\varphi,\omega^*,\omega^\e)$ then the equivalence class of $(a,b,c,\lambda)$ under $\!\fz2s4$ consists of all   elements of ${\F^\times}^4$ which are feasible for $(\mu,\varphi,\omega^*,\omega^\e)$.
\end{enumerate}
\end{thm}
\begin{proof}
(i): Since $\F$ is algebraically closed the statement (i) is immediate from Theorem \ref{thm:universal_class}.

(ii): Since $(a,b,c,\lambda)$ is feasible for $(\mu,\varphi,\omega^*,\omega^\e)$ we may substitute 
Definition \ref{defn:feasible}(i)--(iv) 
into the polynomials given in Theorem \ref{thm:universal_class}(i)--(iii). By Theorem \ref{thm:universal_class} one may factor these polynomials into linear factors to obtain all elements which are feasible for $(\mu,\varphi,\omega^*,\omega^\e)$. By Table \ref{t:z2s4} they are all elements of the equivalence class of $(a,b,c,\lambda)$ under $\!\fz2s4$.
\end{proof}

\section{Proof for Theorem \ref{thm:=d}}
\label{s:=d}

While using the quintuple $(a,b,c,\lambda,\delta)$ to represent an element of ${\F^\times}^4\times \F$ we simply write 
\begin{align*}
w_i=m_i+O_\lambda^\delta(a,b,c)
\qquad 
\hbox{for all $i\in \N$}.
\end{align*}
Recall from Section \ref{s:intro} that $\{w_i\}_{i=0}^{\dbar-1}$ is a basis for $W_\lambda^\delta(a,b,c)$.

\begin{lem}
\label{lem:W_action}
For any $(a,b,c,\lambda,\delta)\in {\F^\times}^4\times \F$ the following statements hold on $W_\lambda^\delta(a,b,c)$: 
\begin{enumerate}
\item The actions of $A$ and $B$ on $W_\lambda^\delta(a,b,c)$ are as follows:
\begin{align*}
(A-\theta_i)w_i &= w_{i+1} \quad (0\leq i\leq \dbar-2),
\qquad 
(A-\theta_{\dbar-1}) w_{\dbar-1}=\delta w_0,
\\
(B-\theta_i^*)w_i &= \varphi_i w_{i-1} \quad (1\leq i\leq \dbar-1),
\qquad 
(B-\theta_0^*)w_0=0.
\end{align*}

\item The elements $\alpha,\beta,\gamma$ act on $W_\lambda^\delta(a,b,c)$ as scalar multiplication by {\rm (\ref{omega})--(\ref{omegae})}, respectively.

\item The element 
$\prod\limits_{i=0}^{\dbar-1}(A-\theta_i)$ acts on $W_\lambda^\delta(a,b,c)$ as scalar multiplication by $\delta$. 
\end{enumerate}
\end{lem}
\begin{proof}
(i): Since $\delta m_0-m_{\dbar}\in O_\lambda^\delta(a,b,c)$ the vector $w_{\dbar}=\delta w_0$. Combined with Theorem \ref{thm:Verma}(i) the statement (i) follows.

(ii): Immediate from Theorem \ref{thm:Verma}(ii).

(iii): Using Lemma \ref{lem:W_action}(i) yields that 
$
\prod_{i=0}^{\dbar-1}(A-\theta_i) w_j
=\delta w_j$ 
for all $j=0,1,\ldots,\dbar-1$. 
Since $W_\lambda^\delta(a,b,c)$ is spanned by $\{w_i\}_{i=0}^{\dbar-1}$ the statement (iii) follows.
\end{proof}

\begin{prop}
\label{prop:W_universal}
Let $(a,b,c,\lambda,\delta)\in {\F^\times}^4\times \F$. 
For any $\triangle_q$-module $V$ and $v\in V$, there exists a $\triangle_q$-module homomorphism $W_\lambda^\delta(a,b,c)\to V$ that sends $w_0$ to $v$ if and only if the following conditions hold:
\begin{enumerate}
\item There exists a $\triangle_q$-module homomorphism 
$
M_\lambda (a,b,c)\to V
$
that sends $m_0$ to $v$.

\item $\prod\limits_{i=0}^{\dbar-1} (A-\theta_i)v=\delta v$.
\end{enumerate}
\end{prop}
\begin{proof}
($\Rightarrow$): Composing the 
$\triangle_q$-module homomorphism $W_\lambda^\delta(a,b,c)\to V$ with the canonical map $M_\lambda(a,b,c)\to W_\lambda^\delta(a,b,c)$ yields the condition (i). The condition (ii) is immediate from Lemma \ref{lem:W_action}(iii).

($\Leftarrow$): By (ii) the vector 
$$
\delta m_0-\prod_{i=0}^{\dbar-1} (A-\theta_i)m_0
$$
lies in the kernel of the $\triangle_q$-module homomorphism $M_\lambda (a,b,c)\to V$ described in (i). By Theorem \ref{thm:Verma}(i) the above subtrahend is equal to $m_{\dbar}$. 
Since the $\triangle_q$-submodule $O_\lambda^{\delta}(a,b,c)$ of $M_\lambda(a,b,c)$ is generated by $\delta m_0-m_{\dbar}$, this induces the desired $\triangle_q$-module homomorphism. The proposition follows.
\end{proof}

For each $n\in \N$ there exists a unique polynomial $T_n(x)\in \F[x]$ such that 
$$
T_n(x+x^{-1})=x^n+x^{-n}.
$$
Note that $T_n(x)$ is of degree $n$ for each integer $n\geq 1$.

\begin{lem}
\label{lem1:Td}
$T_{\dbar}(x)=
\prod\limits_{i=0}^{\dbar-1}(x-\mu q^{2i}-\mu^{-1} q^{-2i})
+\mu^{\dbar}+\mu^{-\dbar}$
for any scalar $\mu\in \F^{\times}$.
\end{lem}
\begin{proof}
Immediate from the construction of $T_{\dbar}(x)$.
\end{proof}

\begin{thm}
[Theorem 3.2, \cite{centerAW:2016}]
\label{thm:Td}
The elements $T_{\dbar}(A)$, $T_{\dbar}(B)$, $T_{\dbar}(C)$ are central in $\triangle_q$.
\end{thm}

\begin{lem}
%[Theorem 3.2, \cite{centerAW:2016}]
\label{lem3:Td}
For any scalar $\mu\in \F^{\times}$ 
the elements 
$$
\prod\limits_{i=0}^{\dbar-1}(A-\mu q^{2i}-\mu^{-1} q^{-2i}),
\quad 
\prod\limits_{i=0}^{\dbar-1}(B-\mu q^{2i}-\mu^{-1} q^{-2i}),
\quad \prod\limits_{i=0}^{\dbar-1}(C-\mu q^{2i}-\mu^{-1} q^{-2i})
$$ 
are central in $\triangle_q$.
\end{lem}
\begin{proof}
Immediate from Lemma \ref{lem1:Td} and Theorem \ref{thm:Td}.
\end{proof}

\begin{proof}[Proof of Theorem \ref{thm:=d}]
Let $\mu$ denote a marginal weight of the $\triangle_q$-module $V$. By Lemma \ref{lem:irr&feasible} there are a nonzero vector $v$ of $V$ and three scalars $\varphi,\omega^*,\omega^\e\in \F$ satisfying the equations (\ref{f&u-1})--(\ref{f&u-4}).

By Theorem \ref{thm:feasible}(i) there exists an element $(a,b,c,\lambda)\in {\F^\times}^4$ which is feasible for $(\mu,\varphi,\omega^*,\omega^\e)$. 
By Theorem \ref{thm:feasible&universal} there exists a unique $\triangle_q$-module homomorphism 
$$
M_\lambda(a,b,c)\to V
$$ 
that sends $m_0$ to $v$. 
By Lemma \ref{lem3:Td} the element 
$
\prod_{i=0}^{\dbar-1}
(A-\theta_i)
$
is central in $\triangle_q$. 
Combined with Schur's lemma there exists a scalar $\delta\in \F$ such that $
\prod_{i=0}^{\dbar-1}
(A-\theta_i)
$ acts on $V$ as scalar multiplication by $\delta$. It follows from Proposition \ref{prop:W_universal} that there exists a unique $\triangle_q$-module homomorphism 
\begin{gather}
\label{W->V}
W_\lambda^\delta (a,b,c)\to V
\end{gather}
that sends $w_0$ to $v$. Since 
$W_\lambda^\delta (a,b,c)$ and $V$ have the same dimension $\dbar$ and by the irreducibility of the $\triangle_q$-module $V$, the map (\ref{W->V}) is a $\triangle_q$-module isomorphism. The result follows.
\end{proof}

\section{The representation-theoretical characterizations of $\!\fz2s4$ and $\!\z2s4$}
\label{s:fz2s4&z2s4}

Starting from this section, we add a few more notational agreements. 
While using the quadruple $(\bar{a},\bar{b},\bar{c},\bar{\lambda})$ to represent an element of ${\F^\times}^4$, the notation $\{\bar{m}_i\}_{i\in \N}$ denotes the basis for $M_{\bar{\lambda}}(\bar{a},\bar{b},\bar{c})$ mentioned in Theorem \ref{thm:Verma}(i) and $\{\bar{\theta}_i\}_{i\in \N}$, $\{\bar{\theta}_i^*\}_{i\in \N}$, $\{\bar{\varphi}_i\}_{i\in \N}$ stand for the accompanying parameters (\ref{thetai})--(\ref{varphii}). While using the quintuple $(\bar{a},\bar{b},\bar{c},\bar{\lambda},\bar{\delta})$ to represent an element of ${\F^\times}^4\times \F$ we write 
\begin{align*}
\bar{w}_i=\bar{m}_i+O_{\bar{\lambda}}^{\bar{\delta}}(\bar{a},\bar{b},\bar{c})
\qquad 
\hbox{for all $i\in \N$}.
\end{align*}

Recall the equivalence relation $\!\fz2s4$ on ${\F^\times}^4$ from Definition \ref{defn:fz2s4&z2s4}(i).

\begin{thm}
\label{thm:fz2s4&M}
For any $(a,b,c,\lambda),(\bar{a},\bar{b},\bar{c},\bar{\lambda})\in {\F^\times}^4$ the following conditions are equivalent:
\begin{enumerate}
\item $(a,b,c,\lambda)\fz2s4\! (\bar{a},\bar{b},\bar{c},\bar{\lambda})$.

\item There exists a $\triangle_q$-module homomorphism
$
M_\lambda(a,b,c)\to M_{\bar{\lambda}}(\bar{a},\bar{b},\bar{c})
$ 
that maps $m_0$ to $\bar{m}_0$.

\item There exists a $\triangle_q$-module homomorphism
$
M_{\bar{\lambda}}(\bar{a},\bar{b},\bar{c}) \to M_\lambda(a,b,c)
$ 
that maps $\bar{m}_0$ to $m_0$.

\item There exists a $\triangle_q$-module isomorphism
$
M_\lambda(a,b,c)\to M_{\bar{\lambda}}(\bar{a},\bar{b},\bar{c})
$ 
that maps $m_0$ to $\bar{m}_0$.

\item There exists a $\triangle_q$-module isomorphism
$
M_{\bar{\lambda}}(\bar{a},\bar{b},\bar{c}) 
\to M_\lambda(a,b,c)
$ 
that maps $\bar{m}_0$ to $m_0$.
\end{enumerate}
\end{thm}
\begin{proof}
(ii), (iii) $\Leftrightarrow$ (iv), (v): Trivial.

(i) $\Rightarrow$ (ii), (iii): By Theorem \ref{thm:feasible}(ii) the condition (i) implies that $(a,b,c,\lambda)$ and $(\bar{a},\bar{b},\bar{c},\bar{\lambda})$ are feasible for the same element of $\F^\times \times \F^3$. 
Combined with Theorem \ref{thm:feasible&universal} the conditions (ii) and (iii) follow.

(ii) $\Rightarrow$ (i): 
By Lemma \ref{lem:Verma_marginalweight}(i) the vector $m_0$ is a marginal weight vector of $M_\lambda(a,b,c)$ with weight $b\lambda^{-1}$.
By Lemma \ref{lem:Verma_marginalweight}(ii) and since $\bar{m}_0$ and $\bar{m}_1$ are linearly independent, the condition (ii) implies $\bar{b}\bar{\lambda}^{-1}=b\lambda^{-1}$. 
Applying Theorem \ref{thm:feasible&universal} yields that $(a,b,c,\lambda)$ and $(\bar{a},\bar{b},\bar{c},\bar{\lambda})$ are feasible for the same element of $\F^\times \times \F^3$. Combined with Theorem \ref{thm:feasible}(ii) the condition (i) follows.

(iii) $\Rightarrow$ (i): Similar to the proof of (ii) $\Rightarrow$ (i).
\end{proof}

\begin{lem}
\label{lem:fz2s4&W}
For any $(a,b,c,\lambda,\delta), (\bar{a},\bar{b},\bar{c},\bar{\lambda},\bar{\delta})\in {\F^\times}^4\times \F$, 
if there exists a $\triangle_q$-module homomorphism 
$
W_\lambda^\delta(a,b,c)
\to
W_{\bar{\lambda}}^{\bar{\delta}}(\bar{a},\bar{b},\bar{c}) 
$
that sends $w_0$ to $\bar{w}_0$ then 
$
(a,b,c,\lambda)\fz2s4\! (\bar{a},\bar{b},\bar{c},\bar{\lambda})$.
\end{lem}
\begin{proof}
Composed with the canonical map $M_\lambda(a,b,c)\to W_\lambda^\delta(a,b,c)$ we obtain the $\triangle_q$-module homomorphism 
$
M_\lambda(a,b,c)
\to
W_{\bar{\lambda}}^{\bar{\delta}}(\bar{a},\bar{b},\bar{c}) 
$
that sends $m_0$ to $\bar{w}_0$.  
By Lemma \ref{lem:Verma_marginalweight}(i) the vector $m_0$ is a marginal weight vector of $M_\lambda(a,b,c)$ with weight $b\lambda^{-1}$.
Since $\bar{w}_0$ and $\bar{w}_1$ are linearly independent and by Lemma \ref{lem:Verma_marginalweight}(ii) this forces that $\bar{b}\bar{\lambda}^{-1}=b\lambda^{-1}$.
Applying Theorem \ref{thm:feasible&universal} yields that $(a,b,c,\lambda)$ and $(\bar{a},\bar{b},\bar{c},\bar{\lambda})$ are feasible for the same element of $\F^\times \times \F^3$. Combined with Theorem \ref{thm:feasible}(ii) the lemma follows.
\end{proof}

\begin{prop}
\label{prop:delta&W}
For any $(a,b,c,\lambda,\delta), (\bar{a},\bar{b},\bar{c},\bar{\lambda},\bar{\delta})\in {\F^\times}^4\times \F$ and any nonzero $\bar{w}\in W_{\bar{\lambda}}^{\bar{\delta}}(\bar{a},\bar{b},\bar{c})$, 
there exists a $\triangle_q$-module homomorphism
$
W_\lambda^\delta (a,b,c)
\to
W_{\bar{\lambda}}^{\bar{\delta}}(\bar{a},\bar{b},\bar{c})
$
that sends $w_0$ to $\bar{w}$ if and only if the following conditions hold:
\begin{enumerate}
\item There exists a $\triangle_q$-module homomorphism 
$
M_\lambda(a,b,c)
\to
W_{\bar{\lambda}}^{\bar{\delta}}(\bar{a},\bar{b},\bar{c})
$
that maps $m_0$ to $\bar{w}$. 

\item The equation {\rm (\ref{delta&bardelta})} holds.
\end{enumerate} 
\end{prop}
\begin{proof}
By Proposition \ref{prop:W_universal} it suffices to show that
\begin{gather}
\label{barthetai_rho(m0)}
\prod_{i=0}^{\dbar-1}(A-\theta_i)\bar{w}
=
\delta\bar{w}
\end{gather}
holds if and only if (\ref{delta&bardelta}) holds.  
To see this we apply the following two equations: 
It follows from Lemma \ref{lem:W_action}(iii) that 
\begin{gather}
\label{thetai_rho(m0)}
\prod_{i=0}^{\dbar-1}(A-\bar{\theta}_i)\bar{w}
=\bar{\delta} \bar{w}.
\end{gather}
It follows from Lemma \ref{lem1:Td} that 
\begin{gather}
\label{thetai&barthetai}
\prod_{i=0}^{\dbar-1}(x-\theta_i)
+a^{\dbar}\lambda^{-\dbar}
+a^{-\dbar}\lambda^{\dbar}
=
\prod_{i=0}^{\dbar-1}(x-\bar{\theta}_i)
+\bar{a}^{\dbar}\bar{\lambda}^{-\dbar}
+\bar{a}^{-\dbar}\bar{\lambda}^{\dbar}.
\end{gather}

(\ref{barthetai_rho(m0)}) $\Rightarrow$ (\ref{delta&bardelta}): We replace $x$ by $A$ and then apply $\bar{w}$ to either side of (\ref{thetai&barthetai}). Simplifying the resulting equation by using (\ref{barthetai_rho(m0)}) and (\ref{thetai_rho(m0)}) it follows that 
$$
(\delta+a^{\dbar}\lambda^{-\dbar}
+a^{-\dbar}\lambda^{\dbar})\bar{w}
=
(\bar{\delta}+\bar{a}^{\dbar}\bar{\lambda}^{-\dbar}
+\bar{a}^{-\dbar}\bar{\lambda}^{\dbar})\bar{w}.
$$
Since $\bar{w}\not=0$ the equation (\ref{delta&bardelta}) follows.

(\ref{delta&bardelta}) $\Rightarrow$ (\ref{barthetai_rho(m0)}): 
Simplifying (\ref{thetai&barthetai}) by using (\ref{delta&bardelta}) yields that  
$$
\prod_{i=0}^{\dbar-1}(x-\theta_i)-\delta
=
\prod_{i=0}^{\dbar-1}(x-\bar{\theta}_i)-\bar{\delta}.
$$ 
We replace $x$ by $A$ and apply $\bar{w}$ to either side of the above equation. Since the right-hand side of the resulting equation is zero by  (\ref{thetai_rho(m0)}), the equation (\ref{barthetai_rho(m0)}) follows.
\end{proof}

Recall the equivalence relation $\!\z2s4$ on ${\F^\times}^4\times \F$ from Definition \ref{defn:fz2s4&z2s4}(ii).

\begin{thm}
\label{thm:z2s4&W}
For any $(a,b,c,\lambda,\delta), (\bar{a},\bar{b},\bar{c},\bar{\lambda},\bar{\delta})\in {\F^\times}^4\times \F$ the following conditions are equivalent:
\begin{enumerate}
\item 
$(a,b,c,\lambda,\delta)
\z2s4\!
(\bar{a},\bar{b},\bar{c},\bar{\lambda},\bar{\delta})$.

\item 
There exists a %unique 
$\triangle_q$-module homomorphism
$
W_\lambda^{\delta}(a,b,c)
\to 
W_{\bar{\lambda}}^{\bar{\delta}}(\bar{a},\bar{b},\bar{c})
$ 
that maps $w_0$ to $\bar{w}_0$.

\item 
There exists a 
$\triangle_q$-module homomorphism 
$
W_{\bar{\lambda}}^{\bar{\delta}}(\bar{a},\bar{b},\bar{c})
\to 
W_\lambda^{\delta}(a,b,c)
$ 
that maps $\bar{w}_0$ to $w_0$.

\item 
There exists a $\triangle_q$-module isomorphism
$
W_\lambda^{\delta}(a,b,c)
\to 
W_{\bar{\lambda}}^{\bar{\delta}}(\bar{a},\bar{b},\bar{c})
$ 
that maps $w_0$ to $\bar{w}_0$.

\item 
There exists a $\triangle_q$-module isomorphism 
$
W_{\bar{\lambda}}^{\bar{\delta}}(\bar{a},\bar{b},\bar{c})
\to 
W_\lambda^{\delta}(a,b,c)
$ 
that maps $\bar{w}_0$ to $w_0$.
\end{enumerate}
\end{thm}
\begin{proof}
(ii), (iii) $\Leftrightarrow$ (iv), (v): Trivial.

(ii), (iii) $\Rightarrow$ (i): 
Immediate from Lemma \ref{lem:fz2s4&W} and Proposition \ref{prop:delta&W}.

(i) $\Rightarrow$ (ii): Since $(a,b,c,\lambda)\fz2s4 \! (\bar{a},\bar{b},\bar{c},\bar{\lambda})$ by Definition \ref{defn:fz2s4&z2s4}(ii), the map mentioned in Theorem \ref{thm:fz2s4&M}(ii) exists. Composed with the canonical map 
$
M_{\bar{\lambda}}(\bar{a},\bar{b},\bar{c})
\to W_{\bar{\lambda}}^{\bar{\delta}}(\bar{a},\bar{b},\bar{c})
$ we obtain the $\triangle_q$-module homomorphism $
M_\lambda(a,b,c)
\to 
W_{\bar{\lambda}}^{\bar{\delta}}(\bar{a},\bar{b},\bar{c})
$
that maps $m_0$ to $\bar{w}_0$. 
Since (\ref{delta&bardelta}) holds by Definition \ref{defn:fz2s4&z2s4}(ii), the condition (ii) is now immediate from Proposition \ref{prop:delta&W}.

(i) $\Rightarrow$ (iii): Similar to the proof of (i) $\Rightarrow$ (ii).
\end{proof}

\section{Proof for Theorem \ref{thm2:irr}}
\label{s:irr}

For convenience we always assume that $(a,b,c,\lambda,\delta)$ denotes a fixed element of ${\F^\times}^4\times \F$ throughout the rest of this paper.

\begin{lem}
\label{lem:necessary}
If the $\triangle_q$-module $W_\lambda^\delta(a,b,c)$ is irreducible then $\delta\not=0$ or $\lambda^2\not\in\{q^{2i}\,|\,i=0,1,\ldots,\dbar-2\}$.
\end{lem}
\begin{proof}
Suppose on the contrary that $\delta=0$ and $\lambda^2=q^{2(i-1)}$ for some integer $i$ with $1\leq i\leq \dbar-1$. 
Let $W$ be the subspace of $W_\lambda^\delta(a,b,c)$ spanned by $\{w_h\}_{h=i}^{\dbar-1}$.
By Lemma \ref{lem:W_action}(ii) the elements $\alpha,\beta,\gamma$ act on $W$ as scalar multiplication.
Since $\delta=0$ and $\varphi_i=0$  by (\ref{varphii}), it follows from Lemma \ref{lem:W_action}(i) that $W$ is invariant under $A$ and $B$.
Hence $W$ is a proper $\triangle_q$-submodule of $W_\lambda^\delta(a,b,c)$ by Proposition \ref{prop:UAW_presentation}, a contradiction. 
\end{proof}

After Theorem \ref{thm:z2s4&W} and Lemma \ref{lem:necessary}, Theorem \ref{thm2:irr}(ii) is already an obvious necessary condition for the irreducibility of the $\triangle_q$-module $W_\lambda^\delta(a,b,c)$.

\begin{proof}[Proof of the implication {\rm (i) $\Rightarrow$ (ii)} of Theorem \ref{thm2:irr}]
In view of Theorem \ref{thm:z2s4&W} the $\triangle_q$-module $W_\lambda^{\delta}(a,b,c)$ is isomorphic to $W_{\bar{\lambda}}^{\bar{\delta}}(\bar{a},\bar{b},\bar{c})$ for any $(\bar{a},\bar{b},\bar{c},\bar{\lambda},\bar{\delta})\z2s4\! (a,b,c,\lambda,\delta)$. 
Combined with Lemma \ref{lem:necessary} the implication (i) $\Rightarrow$ (ii) follows.
\end{proof}

Theorem \ref{thm2:irr}(ii) can be expanded as follows by using Table \ref{t:z2s4}:

\begin{lem}
\label{lem:irr2}
{\rm Theorem \ref{thm2:irr}(ii)} holds if and only if each of the following conditions holds:
\begin{enumerate}
\item  
$\delta\not=0$ or 
$
\lambda^2,
a^{-1}b^{-1}c^{-1}\lambda q^{-1},
a^{-1}b^{-1}c\lambda q^{-1}
\not\in \{q^{2i}\,|\,i=0,1,\ldots,\dbar-2\}.
$

\item 
$
\delta\not=(a^{\dbar}-a^{-\dbar})(\lambda^{\dbar}-\lambda^{-\dbar})$ 
or 
$$
\lambda^2,
ab^{-1}c^{-1}\lambda q^{-1},
ab^{-1}c\lambda q^{-1}
\not\in \{q^{2i}\,|\,i=0,1,\ldots,\dbar-2\}.
$$

\item 
$
\delta+a^{\dbar}\lambda^{-\dbar}+a^{-\dbar}\lambda^{\dbar}\not=(b^{\dbar}c^{\dbar}+b^{-\dbar}c^{-\dbar})q^{\dbar}$
or 
$$
ab^{-1}c^{-1}\lambda q^{-1},
a^{-1}b^{-1}c^{-1}\lambda q^{-1},
b^{-2} q^{-2}
\not\in \{q^{2i}\,|\,i=0,1,\ldots,\dbar-2\}.
$$

\item 
$
\delta+a^{\dbar}\lambda^{-\dbar}+a^{-\dbar}\lambda^{\dbar}\not=(b^{\dbar}c^{-\dbar}+b^{-\dbar}c^{\dbar})q^{\dbar}$
or 
$$
ab^{-1}c\lambda q^{-1},
b^{-2}q^{-2},
a^{-1}b^{-1}c\lambda q^{-1}
\not\in \{q^{2i}\,|\,i=0,1,\ldots,\dbar-2\}.
$$
\end{enumerate}
\end{lem}

Recall the algebra automorphism $\vee$ of $\triangle_q$ from Lemma \ref{lem:vee}. We are going to investigate the marginal weights and the marginal weight vectors of $W_\lambda^\delta(a,b,c)^\vee$.
Define $\nu$ to be a scalar in $\F^\times$ satisfying 
\begin{gather}
\label{W:nu}
\nu^{\dbar} +\nu^{-\dbar}=\delta+a^{\dbar} \lambda^{-\dbar}+a^{-\dbar} \lambda^{\dbar}.
\end{gather}
Since $\F$ is algebraically closed the existence of $\nu$ follows. Set 
$$
\vartheta_i=
\nu^{-1} q^{2i} 
+
\nu q^{-2i}
\qquad \hbox{for all integers $i$}. 
$$

\begin{lem}
\label{lem2:minpoly_W}
The characteristic polynomial of $A$ on $W_\lambda^\delta(a,b,c)$ is equal to 
\begin{gather}\label{Acharapoly_W}
\prod\limits_{i=0}^{\dbar-1}(x-\vartheta_i)=\prod_{i=0}^{\dbar-1}(x-\theta_i)-\delta.
\end{gather}
\end{lem}
\begin{proof}
By Lemma \ref{lem:W_action}(iii) the characteristic polynomial of $A$ on $W_\lambda^\delta(a,b,c)$ is equal to 
the right-hand side of (\ref{Acharapoly_W}). 
It follows from Lemma \ref{lem1:Td} that 
$$
\prod_{i=0}^{\dbar-1}(x-\theta_i)+a^{\dbar}\lambda^{-\dbar}+a^{-\dbar}\lambda^{\dbar}=\prod_{i=0}^{\dbar-1}(x-\vartheta_i)+\nu^{\dbar}+\nu^{-\dbar}.
$$
Combined with (\ref{W:nu}) the lemma follows.
\end{proof}

For any integer $i$ with $0\leq i\leq \dbar-1$ we define $e_i\in W_\lambda^\delta(a,b,c)$ by 
\begin{gather}\label{ei}
e_i=\sum_{h=1}^{\dbar}
\prod_{j=h}^{\dbar-1}(\vartheta_i-\theta_j) w_{h-1}.
\end{gather}

\begin{lem}
\label{lem:ei}
For any integer $i$ with $0\leq i\leq \dbar-1$ the coefficient of $w_{\dbar-1}$ in $e_i$ with respect to the basis $\{w_h\}_{h=0}^{\dbar-1}$ for $W_\lambda^\delta(a,b,c)$ is equal to one.
\end{lem}
\begin{proof}
Immediate from (\ref{ei}).
\end{proof}

By Lemma \ref{lem2:minpoly_W} the weights of $W_\lambda^\delta(a,b,c)^\vee$ are $\{\nu q^{-2i}\}_{i=0}^{\dbar-1}$ and $\{\nu^{-1} q^{2i}\}_{i=0}^{\dbar-1}$.

\begin{lem}
\label{lem2:weight}
For any integer $i$  with $0\leq i\leq \dbar-1$ the following statements are true:
\begin{enumerate}
\item The weight space of $W_\lambda^\delta(a,b,c)^\vee$ with weight $\nu^{-1} q^{2i}$ is spanned by $e_i$. 

\item The weight space of $W_\lambda^\delta(a,b,c)^\vee$ with weight $\nu q^{-2i}$ is spanned by $e_i$.
\end{enumerate}
\end{lem}
\begin{proof} 
Applying Lemma \ref{lem:W_action}(i) a direct calculation shows that $(A-\vartheta_i)e_i$ is equal to 
$$
\delta w_0
+\sum_{h=1}^{\dbar-1}
\prod_{j=h}^{d-1} (\vartheta_i-\theta_j) w_h
-\sum_{h=1}^{\dbar}
\prod_{j=h-1}^{d-1}(\vartheta_i-\theta_j) w_{h-1}.
$$
Using the change of indices yields that the above vector is equal to the scalar multiple of $w_0$ by 
$$
\delta-\prod_{j=0}^{\dbar-1}(\vartheta_i-\theta_j), 
$$
which is equal to zero by (\ref{Acharapoly_W}). Therefore $e_i$ is a weight vector of $W_\lambda^\delta(a,b,c)^\vee$ with weights $\nu q^{-2i}$ and $\nu^{-1} q^{2i}$. 
By Lemma \ref{lem:W_action}(i) the rank of $A-\vartheta_i$ on $W_\lambda^\delta(a,b,c)$ is at least $\dbar-1$. The lemma now follows from the rank-nullity theorem.
\end{proof}

\begin{lem}
\label{lem1:marginal}
For any integer $i$ with $0\leq i\leq \dbar-1$ the following conditions are equivalent:
\begin{enumerate}
\item $\nu^{-1} q^{2i}$ is a marginal weight of  $W_\lambda^\delta(a,b,c)^\vee$.

\item $(A-\vartheta_{i+1})(A-\vartheta_i)Be_i=0$ on $W_\lambda^\delta (a,b,c)$.

\item $e_i$ is a marginal weight vector of  $W_\lambda^\delta(a,b,c)^\vee$ with weight $\nu^{-1} q^{2i}$.

\item $\nu\in 
\{
a \lambda^{-1} q^{2(i-1)},
a^{-1} \lambda^{-1} q^{2(i-1)},
b c q^{2i-1},
b c^{-1} q^{2i-1} 
\}$.
\end{enumerate}
\end{lem}
\begin{proof} 
(i) $\Leftrightarrow$ (ii): Immediate from Definition \ref{defn:marginal} and Lemma \ref{lem2:weight}(i).

(ii) $\Leftrightarrow$ (iii): Immediate from Definition \ref{defn:marginalvector} and Lemma \ref{lem2:weight}(i).

(ii) $\Leftrightarrow$ (iv):  
%By Lemmas \ref{lem:weight}(ii) and \ref{lem2:weight}(i)  the vector
Applying Lemma \ref{lem:weight}(ii) to $W_\lambda^\delta(a,b,c)^\vee$ yields that 
\begin{gather}
\label{AAei}
(A-\vartheta_{i+1})(A-\vartheta_i)Be_i\in W_\lambda^\delta(a,b,c)^\vee(\nu^{-1}q^{2(i-1)}).
\end{gather}
Combined with Lemma \ref{lem2:weight}(i) the left-hand side of (\ref{AAei}) is a scalar multiple of $e_{i-1}$, where $e_{i-1}$ is interpreted as $e_{\dbar-1}$ if $i=0$.
In view of Lemma \ref{lem:ei} 
the condition (ii) holds if and only if the coefficient of $w_{\dbar-1}$ in the left-hand side of (\ref{AAei}) is equal to zero.
Applying Lemma \ref{lem:W_action}(i) a direct calculation shows that the aforementioned coefficient is equal to $b^{-1}\lambda^{-1} q^{-1}(q-q^{-1})(q^2-q^{-2})$ times 
\begin{align*}
%&b^{-1}\lambda^{-1} q^{-1}(q-q^{-1})(q^2-q^{-2})
%\\
%&\quad \times\; 
(q^{1-i}-bc \nu^{-1} q^i) (q^{1-i}-bc^{-1} \nu^{-1} q^i)
(q^{i-1}-a \lambda \nu q^{1-i}) (q^{i-1}-a^{-1} \lambda \nu q^{1-i}).
\end{align*}
The equivalence of (ii) and (iv) follows. 
\end{proof}

\begin{lem}
\label{lem2:marginal}
For any integer $i$ with $0\leq i\leq \dbar-1$ the following conditions are equivalent:
\begin{enumerate}
\item $\nu q^{-2i}$ is a marginal weight of $W_\lambda^\delta(a,b,c)^\vee$.

\item $(A-\vartheta_{i-1})(A-\vartheta_i)Be_i=0$ on $W_\lambda^\delta (a,b,c)$.

\item $e_i$ is a marginal weight vector of $W_\lambda^\delta(a,b,c)^\vee$ with weight $\nu q^{-2i}$.

\item  $\nu \in 
\{
a \lambda q^{2(i+1)},
a^{-1} \lambda q^{2(i+1)},
b^{-1} c q^{2i+1}, 
b^{-1} c^{-1} q^{2i+1}
\}$.
\end{enumerate}
\end{lem}
\begin{proof}
Similar to the proof of Lemma \ref{lem1:marginal}.
\end{proof}

\begin{lem}
\label{lem3:marginal}
For any integer $i$ with $0\leq i\leq \dbar-1$ the following conditions are equivalent:
\begin{enumerate}
\item $e_i$ is a marginal weight vector of the $\triangle_q$-module $W_\lambda^\delta(a,b,c)^\vee$.

\item $\nu q^{-2i}\in\{
a \lambda^{-1} q^{-2},
a^{-1} \lambda q^{2},
a \lambda q^{2},
a^{-1}\lambda^{-1} q^{-2},
b c q^{-1},
b^{-1} c^{-1} q,
b c^{-1} q^{-1},
b^{-1} c q\}$.

\item $\nu^{-1} q^{2i}\in\{
a \lambda^{-1} q^{-2},
a^{-1} \lambda q^{2},
a \lambda q^{2},
a^{-1}\lambda^{-1} q^{-2},
b c q^{-1},
b^{-1} c^{-1} q,
b c^{-1} q^{-1},
b^{-1} c q\}$.
\end{enumerate}
\end{lem}
\begin{proof}
Combine Lemmas \ref{lem1:marginal} and \ref{lem2:marginal}.
\end{proof}

Let $L_{jk}^{(i)}$ for all integers $i,j,k$ with $0\leq i,j,k\leq \dbar-1$ denote the scalars in $\F$ satisfying 
\begin{gather}
\label{Lijk}
\prod_{h=1}^k (B-\theta_{\dbar-h}^*) e_i
=\sum_{j=0}^{\dbar-1} L_{jk}^{(i)} w_{\dbar-j-1}.
\end{gather}
Let $i$ denote an integer with $0\leq i\leq \dbar-1$.
Using Lemma \ref{lem:W_action}(i) yields that 
$L_{jk}^{(i)}=0$ for all integers $j,k$ with 
$0\leq j<k\leq \dbar-1$. 
In particular 
\begin{gather}
\label{Ld-1d-1w0}
\prod_{h=1}^{\dbar-1} (B-\theta_{\dbar-h}^*) e_i
=L_{\dbar-1,\dbar-1}^{(i)} w_0.
\end{gather}

\begin{thm}
\label{thm1:irr}
The $\triangle_q$-module $W_\lambda^\delta(a,b,c)$ is reducible if and only if there exists an integer $i$ with $0\leq i\leq \dbar-1$ satisfying the following conditions: 
\begin{enumerate}
\item $e_i$ is a marginal weight vector of the $\triangle_q$-module $W_\lambda^\delta(a,b,c)^\vee$.

\item $L_{\dbar-1,\dbar-1}^{(i)}=0$.
\end{enumerate}
\end{thm}
\begin{proof}
($\Rightarrow$): Since the $\triangle_q$-module $W_\lambda^\delta(a,b,c)$ is reducible, there is a proper irreducible $\triangle_q$-submodule $W$ of $W_\lambda^\delta(a,b,c)$. 
By Lemma \ref{lem:<d&marginal} there exists a marginal weight $\mu$ of $W^\vee$. 
By Lemma \ref{lem2:minpoly_W} there exists an integer $i$ with $0\leq i\leq \dbar-1$ such that $\mu=\nu q^{-2i}$ or $\mu=\nu^{-1}q^{2i}$.
Combined with Lemmas \ref{lem1:marginal} and \ref{lem2:marginal} the vector $e_i$ is a marginal weight vector of $W^\vee$. The condition (i) follows.
By (\ref{Ld-1d-1w0}) the vector $L_{\dbar-1,\dbar-1}^{(i)} w_0\in W$.
Since the $\triangle_q$-module $W_\lambda^\delta(a,b,c)$ is generated by $w_0$ the vector $w_0\not\in W$. The condition (ii) follows.

($\Leftarrow$): 
Suppose on the contrary that the $\triangle_q$-module $W_\lambda^\delta(a,b,c)$ is irreducible. 
Applying Lemma \ref{lem3:marginalvector} to $W_\lambda^\delta(a,b,c)^\vee$ the vectors $\{B^h e_i\}_{h=0}^{\dbar-1}$ form a basis for $W_\lambda^\delta(a,b,c)$.
 By (ii) the right-hand side of (\ref{Ld-1d-1w0}) is zero. Then $\{B^h e_i\}_{h=0}^{\dbar-1}$ are linearly dependent, a contradiction.
\end{proof}

\begin{prop}
\label{prop:Ljk}
For any integer $i$ with $0\leq i\leq \dbar-1$ the following statements are true:
\begin{enumerate}
\item Suppose that $\nu  q^{-2i}\in\{a \lambda^{-1} q^{-2}, a^{-1} \lambda q^{2}\}$. Then 
\begin{align*}
L_{jk}^{(i)}=\left\{
\begin{array}{cl}
\displaystyle
\prod_{h=1}^j
\varphi_{\dbar-h}
%(q^{h}-a b c^{-1} \lambda^{-1} q^{-h-1})
%(q^{h}-a b c \lambda^{-1} q^{-h-1})
%(\lambda q^{h+1}-\lambda^{-1} q^{-h-1})
\qquad 
&\hbox{if $j=k$},
\\
0
\qquad 
&\hbox{if $j\not=k$}
\end{array}
\right.
\qquad 
(0\leq j,k\leq \dbar-1).
\end{align*}

\item Suppose that $\nu q^{-2i}\in\{ a\lambda q^{2},
a^{-1} \lambda^{-1} q^{-2}\}$. Then $L_{jk}^{(i)}$ is equal to 
\begin{align*}
\begin{split}
&\prod_{h=1}^k
(q^{h+j-k}-q^{k-h-j})
(a \lambda q^{h+1}-b c q^{-h})
(b^{-1} q^h-a^{-1} c^{-1} \lambda^{-1} q^{-h-1})
\\
&\qquad \times \,
\prod_{h=1}^j 
(\lambda q^{h+1}-\lambda^{-1} q^{-h-1})
\prod_{h=1}^{j-k} 
(a q^{1-h}-a^{-1} q^{h-1}) 
\end{split}
\qquad 
(0\leq j,k\leq \dbar-1).
\end{align*}

\item Suppose that $\nu q^{-2i}\in\{b c q^{-1},b^{-1} c^{-1} q\}$. Then $L_{jk}^{(i)}$ is equal to 
\begin{align*}
\begin{split}
&\prod_{h=1}^k
(q^{h+j-k}-q^{k-h-j})
(b q^{-h}-b^{-1} q^h)
(a \lambda q^{h+1}-b c q^{-h})
\\
&\qquad \times \,
\prod_{h=1}^j
(\lambda^{-1} q^{-h-1}-a^{-1} b^{-1} c^{-1} q^{h})
\prod_{h=1}^{j-k}
(b c\lambda q^{h}-a q^{1-h})
\end{split}
\qquad 
(0\leq j,k\leq \dbar-1).
\end{align*}

\item Suppose that $\nu q^{-2i}\in\{b c^{-1} q^{-1},b^{-1} c q\}$. Then $L_{jk}^{(i)}$ is equal to 
\begin{align*}
\begin{split}
&\prod_{h=1}^k
(q^{h+j-k}-q^{k-h-j})
(b q^{-h}-b^{-1} q^h)
(a \lambda q^{h+1}-b c^{-1} q^{-h})
\\
&\qquad \times \,
\prod_{h=1}^j
(\lambda^{-1} q^{-h-1}-a^{-1} b^{-1} c q^{h})
\prod_{h=1}^{j-k}
(b c^{-1} \lambda q^{h}-a q^{1-h})
\end{split}
\qquad 
(0\leq j,k\leq \dbar-1).
\end{align*}
\end{enumerate}
\end{prop}
\begin{proof}
Fix an integer $i$ with $0\leq i\leq \dbar-1$. 
Recall the coefficients $\{L_{jk}^{(i)}\}_{0\leq j,k\leq \dbar-1}$ from (\ref{Lijk}).
Applying Lemma \ref{lem:W_action}(i) yields that  
\begin{align}
\label{L0k}
L_{0k}^{(i)} &=0
\qquad 
(1\leq k\leq \dbar-1),
\\
\label{Ljk:rec}
L_{jk}^{(i)}&=
\varphi_{\dbar-j} L_{j-1,k-1}^{(i)}
+
(\theta_{\dbar-j-1}^*-\theta_{\dbar-k}^*) L_{j,k-1}^{(i)}
\qquad 
(1\leq j,k \leq \dbar-1).
\end{align}
By (\ref{ei}) we have
\begin{gather}
\label{Lj0}
L_{j0}^{(i)}=\prod_{h=1}^j(\vartheta_i-\theta_{\dbar-h})
\qquad 
(0\leq j\leq \dbar-1).
\end{gather}
The coefficients $\{L_{jk}^{(i)}\}_{0\leq j,k\leq \dbar-1}$ are uniquely determined by the recurrence (\ref{Ljk:rec}) and the initial conditions (\ref{L0k}) and (\ref{Lj0}).  To see (i)--(iv) it is routine but tedious to verify that the given formulae satisfy (\ref{L0k})--(\ref{Lj0}).
\end{proof}

\begin{proof}[Proof of the implication {\rm (ii) $\Rightarrow$ (i)} of Theorem \ref{thm2:irr}] 
Suppose on the contrary that Theorem \ref{thm2:irr}(i) fails. By Theorem \ref{thm1:irr} there is an integer $i$ with $0\leq i\leq \dbar-1$ such that Theorem \ref{thm1:irr}(i), (ii) hold. By Lemma \ref{lem3:marginal} there are four possible situations: 
(a) $\nu q^{-2i}\in\{a\lambda^{-1}q^{-2},a^{-1}\lambda q^2\}$; 
(b) $\nu q^{-2i}\in\{a\lambda q^2,a^{-1}\lambda^{-1} q^{-2}\}$;
(c) $\nu q^{-2i}\in\{bcq^{-1},b^{-1}c^{-1}q\}$;
(d) $\nu q^{-2i}\in\{bc^{-1}q^{-1},b^{-1}cq\}$.

(a): Since the order of $q^2$ is $\dbar$ the scalar
$\nu^{\dbar}\in \{a^{\dbar}\lambda^{-\dbar},a^{-\dbar}\lambda^{\dbar}\}$. It follows from (\ref{W:nu}) that $\delta=0$. By Lemma \ref{lem:irr2}(i) this forces that
$$
\lambda^2,a^{-1}b^{-1}c^{-1}\lambda q^{-1},a^{-1}b^{-1}c\lambda q^{-1}\not\in \{q^{2i}\,|\,i=0,1,\ldots,\dbar-2\}.
$$
Hence $\varphi_i\not=0$ for all $i=1,2,\ldots,\dbar-1$ by (\ref{varphii}).
Combined with Proposition \ref{prop:Ljk}(i) this yields that 
$L_{\dbar-1,\dbar-1}^{(i)}\not=0$, a contradiction to Theorem \ref{thm1:irr}(ii).

(b): Since the order of $q^2$ is $\dbar$ the scalar
$\nu^{\dbar}\in \{a^{\dbar}\lambda^{\dbar},a^{-\dbar}\lambda^{-\dbar}\}$. It follows from (\ref{W:nu}) that $\delta=(a^{\dbar}-a^{-\dbar})(\lambda^{\dbar}-\lambda^{-\dbar})$. By Lemma \ref{lem:irr2}(ii) this forces that
$$
\lambda^2, 
ab^{-1}c^{-1}\lambda q^{-1}, 
ab^{-1}c\lambda q^{-1}
\not\in \{q^{2i}\,|\,i=0,1,\ldots,\dbar-2\}.
$$
%Combined with Proposition \ref{prop:Ljk}(ii) this yields that 
Then $L_{\dbar-1,\dbar-1}^{(i)}\not=0$ by Proposition \ref{prop:Ljk}(ii), a contradiction to Theorem \ref{thm1:irr}(ii).

(c): Since the order of $q^2$ is $\dbar$ the scalar
$\nu^{\dbar}\in \{b^{\dbar}c^{\dbar}q^{\dbar}, b^{-\dbar}c^{-\dbar}q^{\dbar}\}$. It follows from (\ref{W:nu}) that $\delta+a^{\dbar}\lambda^{-\dbar}+a^{-\dbar}\lambda^{\dbar}=(b^{\dbar}c^{\dbar}+b^{-\dbar}c^{-\dbar})q^{\dbar}$. By Lemma \ref{lem:irr2}(iii) this forces that
$$
ab^{-1}c^{-1}\lambda q^{-1}, 
a^{-1}b^{-1}c^{-1}\lambda q^{-1},
b^{-2}q^{-2}
\not\in \{q^{2i}\,|\,i=0,1,\ldots,\dbar-2\}.
$$
%Combined with Proposition \ref{prop:Ljk}(iii) this yields that 
Then 
$L_{\dbar-1,\dbar-1}^{(i)}\not=0$ by Proposition \ref{prop:Ljk}(iii), a contradiction to Theorem \ref{thm1:irr}(ii).

(d): Since the order of $q^2$ is $\dbar$ the scalar
$\nu^{\dbar}\in \{b^{\dbar}c^{-\dbar}q^{\dbar},b^{-\dbar}c^{\dbar}q^{\dbar}\}$. It follows from (\ref{W:nu}) that $\delta+a^{\dbar}\lambda^{-\dbar}+a^{-\dbar}\lambda^{\dbar}=(b^{\dbar}c^{-\dbar}+b^{-\dbar}c^{\dbar})q^{\dbar}$. By Lemma \ref{lem:irr2}(iv) this forces that
$$
ab^{-1}c\lambda q^{-1}, 
b^{-2}q^{-2},
a^{-1}b^{-1}c\lambda q^{-1}
\not\in \{q^{2i}\,|\,i=0,1,\ldots,\dbar-2\}.
$$
%Combined with Proposition \ref{prop:Ljk}(iv) this yields that 
Then 
$L_{\dbar-1,\dbar-1}^{(i)}\not=0$ 
by Proposition \ref{prop:Ljk}(iv), a contradiction to Theorem \ref{thm1:irr}(ii).
\end{proof}

\section{The binary relation $\sim$ and the marginal weight vectors of $W_\lambda^\delta(a,b,c)$}
\label{s:sim}

\begin{lem}
\label{lem:BeigenW}
The characteristic polynomial of $B$ on $W_\lambda^\delta(a,b,c)$ is equal to 
$$
\prod_{i=0}^{\dbar-1} 
(x-\theta_i^*).
$$
\end{lem}
\begin{proof}
Immediate from Lemma \ref{lem:W_action}(i).
\end{proof}

For any integers $i,j$ with $0\leq i\leq j\leq \dbar-1$ we define $w_{ij}\in W_{\lambda}^\delta(a,b,c)$ by 
\begin{gather}
\label{wij}
w_{ij}=\sum_{h=0}^{j-i}
\left(\prod_{k=h}^{j-i-1}\varphi_{i+k+1}\right)
\left(\prod_{k=0}^{h-1}(\theta_j^*-\theta_{i+k}^*)\right)
w_{i+h}.
\end{gather}
Note that $w_{ii}=w_i$ for all $i=0,1,\ldots,\dbar-1$.  
Applying Lemma \ref{lem:W_action}(i) a straightforward calculation yields the following three lemmas:

\begin{lem}
\label{lem1:Wmarginal}
If there exists an integer $i$ with $0\leq i\leq \dbar-1$ such that $\varphi_i=0$, then 
$$
(B-\theta_j^*)w_{ij}=0
$$ 
for any integer $j$ with $i\leq j\leq \dbar-1$. 
In particular $(B-\theta_i^*)w_{0i}=0$ for any integer $i$ with $0\leq i\leq \dbar-1$.
\end{lem}

\begin{lem}
\label{lem3:Wmarginal}
For each integer $i$ with $1\leq i\leq \dbar-1$ 
the vector $(B-\theta_{i+1}^*)(B-\theta_i^*)A w_{0i}$ is equal to the scalar multiple of $w_{0,i-1}$ by 
\begin{align*}
&a b^{-1}\lambda q (q-q^{-1}) (q^2-q^{-2}) 
(q^i-q^{-i})\cdot \varphi_i
\\
&\quad \times\,
(b q^i-b^{-1} q^{-i})
(q^{-i}-a^{-1}bc^{-1} \lambda^{-1} q^{i-1})
(q^{-i}-a^{-1}bc \lambda^{-1} q^{i-1}).
\end{align*}
\end{lem}

\begin{lem}
\label{lem2:Wmarginal}
The following statements hold on $W_\lambda^\delta(a,b,c)$:
\begin{enumerate}
\item $(B-\theta_{\dbar-1}^*)(B-\theta_0^*)A w_0=(\theta_1^*-\theta_{\dbar-1}^*)w_{01}$.

\item For each integer $i$ with $1\leq i\leq \dbar-2$ the vector $(B-\theta_{i-1}^*)(B-\theta_i^*)A w_{0i}$ is equal to 
\begin{gather*}
\frac{(q-q^{-1})(q^2-q^{-2})}{(q^i-q^{-i})(q^{i+1}-q^{-i-1})}(\theta_i^*-\theta_0^*) w_{0,i+1}.
\end{gather*}

\item $(B-\theta_{\dbar-2}^*)(B-\theta_{\dbar-1}^*)A w_{0,\dbar-1}$ is equal to the scalar multiple of $w_0$ by 
\begin{align*}
&q^{-\frac{\dbar(\dbar-1)}{2}}(q^2-q^{-2})\prod_{i=1}^{\dbar-1}(q^i-q^{-i})
\cdot 
(\theta_0^*-\theta_{\dbar-1}^*)
\\
&\quad\times\,
\left(
\delta (b^{\dbar}\lambda^{-\dbar}-b^{-\dbar}\lambda^{\dbar})-a^{-\dbar}b^{-\dbar}
(\lambda^{2\dbar}-1)
(a^{\dbar}b^{\dbar}c^{\dbar}\lambda^{-\dbar}q^{\dbar}-1)
(a^{\dbar}b^{\dbar}c^{-\dbar}\lambda^{-\dbar}q^{\dbar}-1)
\right)\!.
\end{align*} 
\end{enumerate}
\end{lem}

\begin{lem}
\label{lem1-2:Wmarginal}
Suppose that there are two integers $i,j$ with $0\leq i\leq j\leq \dbar-1$ such that $\varphi_i=0$ and one of the following conditions holds:
\begin{enumerate}
\item $\varphi_{i+1}\varphi_{i+2}\cdots \varphi_j\not=0$.

\item $\theta_j^*\not\in\{\theta_h^*\,|\,h=i,i+1,\ldots,j-1\}$. 
\end{enumerate}
Then $w_{ij}$ is a basis for the $\theta_j^*$-eigenspace of $B$ on the subspace of $W_\lambda^\delta(a,b,c)$ spanned by $\{w_h\}_{h=i}^j$. 
\end{lem}
\begin{proof}
Let $W$ denote the subspace of $W_\lambda^\delta(a,b,c)$ spanned by $\{w_h\}_{h=i}^j$. Since $\varphi_i=0$ it follows from Lemma \ref{lem:W_action}(i) that $W$ is $B$-invariant and the characteristic polynomial of $B$ on $W$ is $\prod_{h=i}^j(x-\theta_h^*)$.

Suppose that (i) holds. Then the rank of $B-\theta_j^*$ on $W$ is equal to $j-i$. By the rank-nullity theorem the $\theta_j^*$-eigenspace of $B$ on $W$ is of dimension one.  
By (\ref{wij}) the coefficient of $w_i$ in $w_{ij}$ is $\varphi_{i+1}\varphi_{i+2}\cdots \varphi_j\not=0$. Hence $w_{ij}$ is nonzero. Combined with Lemma \ref{lem1:Wmarginal} the lemma  is true when (i) holds.

Suppose that (ii) holds. Since $\theta_j^*$ is a simple root of $\prod_{h=i}^j(x-\theta_h^*)$, the $\theta_j^*$-eigenspace of $B$ on $W$ is of dimension one. By (\ref{wij}) the coefficient of $w_j$ in $w_{ij}$ is $\prod_{h=i}^{j-1}(\theta_j^*-\theta_h^*)\not=0$. Hence $w_{ij}$ is nonzero. Combined with Lemma \ref{lem1:Wmarginal} the lemma is true when (ii) holds.
\end{proof}

Recall the binary relation $\sim$ on ${\F^\times}^4\times \F$ from Definition \ref{defn:sim}.

\begin{lem}
\label{lem1-0:existence} 
Suppose that there exists an integer $i$ with $1\leq i\leq \dbar-1$
such that 
$\lambda^2=q^{2(i-1)}$. 
Then $
(\bar{a},\bar{b},\bar{c},\bar{\lambda},\bar{\delta})=(a^{-1},b,c,\lambda^{-1}q^{-2},\delta)$  satisfies following statements:
\begin{enumerate}
\item $\bar{\theta}_h^*=\theta_{i+h}^*$ for all $h\in \N$. 

\item $w_i$ is a marginal weight vector of $W_\lambda^\delta(a,b,c)$ with weight $\bar{b}\bar{\lambda}^{-1}$. 

\item There is a $\triangle_q$-module homomorphism $W_{\bar{\lambda}}^{\bar{\delta}}(\bar{a},\bar{b},\bar{c})\to W_\lambda^\delta(a,b,c)$ that maps 
$\bar{w}_0$ to $w_i$.
%\begin{eqnarray*}
%\bar{w}_0 &\mapsto & w_i.
%\end{eqnarray*}
\end{enumerate}
\end{lem}
\begin{proof}
Since $\lambda^2=q^{2(i-1)}$ and $(\bar{b},\bar{\lambda})=(b,\lambda^{-1}q^{-2})$ it follows that  
\begin{gather}
\label{lem1-0:existence:mu}
\bar{b}\bar{\lambda}^{-1}=b \lambda^{-1} q^{2i}.
\end{gather}
The condition (i) follows. 
Let $\mu$ denote the scalar (\ref{lem1-0:existence:mu}).
Then $(\mu+\mu^{-1},\mu q^2+\mu^{-1} q^{-2})=(\theta_i^*,\theta_{i+1}^*)$. 
Since $\lambda^2=q^{2(i-1)}$ it follows from (\ref{varphii}) that $\varphi_i=0$. 
By Lemma \ref{lem1:Wmarginal} the vector $w_i$ is a weight vector of $W_\lambda^{\delta}(a,b,c)$ with weight $\mu$. 
Let 
$$
\varphi=\frac{\varphi_{i+1}+\theta_i(\theta_i^*-\theta_{i+1}^*)}{q-q^{-1}}.
$$ 
Using Lemma \ref{lem:W_action}(i) a routine calculation yields that 
$
(B-\theta_{i+1}^*)Aw_i=(q-q^{-1})\varphi\cdot w_i.
$ 
The condition (ii) follows.

Under the hypothesis $\lambda^2=q^{2(i-1)}$ the element $(\bar{a},\bar{b},\bar{c},\bar{\lambda})$ is feasible for $(\mu,\varphi,\omega^*,\omega^\e)$ where 
$\omega^*$ and $\omega^{\e}$ are the scalars (\ref{omega*}) and (\ref{omegae}) respectively. 
By Theorem \ref{thm:feasible&universal} there exists a $\triangle_q$-module homomorphism 
$
M_{\bar{\lambda}}(\bar{a},\bar{b},\bar{c})\to W_\lambda^{\delta}(a,b,c)
$
that sends $\bar{m}_0$ to $w_i$. 
Since $(\bar{a}^{\dbar}\bar{\lambda}^{-\dbar},\bar{\delta})=(a^{-\dbar}\lambda^{\dbar},\delta)$ the equation (\ref{delta&bardelta}) holds. 
By Proposition \ref{prop:delta&W} the condition (iii) follows.
\end{proof}

\begin{lem}
\label{lem1-2:existence}
Suppose that the conditions {\rm (a)} and {\rm (b)} of {\rm Definition \ref{defn:sim}(iii)} hold. Then $(\bar{a},\bar{b},\bar{c},\bar{\lambda},\bar{\delta})=(a^{-1},b^{-1},c,\lambda^{-1}q^{-2},\delta)$
satisfies the following statements:
\begin{enumerate}
\item $\bar{\theta}_h^*=\theta_{\dbar-h-1}^*$ for all  $h=0,1,\ldots, \dbar-1$.

\item $w_{0,\dbar-1}$ is a marginal weight vector of $W_\lambda^\delta(a,b,c)$ with weight $\bar{b}\bar{\lambda}^{-1}$.

\item There is a $\triangle_q$-module homomorphism $W_{\bar{\lambda}}^{\bar{\delta}}(\bar{a},\bar{b},\bar{c})\to W_\lambda^\delta(a,b,c)$ that maps $\bar{w}_0$ to $w_{0,\dbar-1}$.
%\begin{eqnarray*}
%\bar{w}_0 &\mapsto & w_{0,\dbar-1}.
%\end{eqnarray*}
\end{enumerate}
\end{lem}
\begin{proof}
Since $(\bar{b},\bar{\lambda})=(b^{-1},\lambda^{-1}q^{-2})$ it follows that
\begin{gather}
\label{lem1-2:existence:mu}
\bar{b}\bar{\lambda}^{-1}=b^{-1}\lambda q^2.
\end{gather}
The condition (i) follows. 
Let $\mu$ denote the scalar (\ref{lem1-2:existence:mu}). Then $(\mu+\mu^{-1},\mu q^2+\mu^{-1} q^{-2})=(\theta_{\dbar-1}^*,\theta_{\dbar-2}^*)$. 
Observe that the condition (a) of Definition \ref{defn:sim}(iii) is equivalent to 
\begin{gather}
\label{condition:a}
\theta_{\dbar-1}^*\not=\theta_i^*
\qquad 
\hbox{for all $i=0,1,\ldots,\dbar-2$}.
\end{gather}
By Lemma \ref{lem1-2:Wmarginal} the vector $w_{0,\dbar-1}$ is a basis for $W_\lambda^\delta(a,b,c)(\mu)$.   
Under the condition (b) of Definition \ref{defn:sim}(iii), Lemma \ref{lem2:Wmarginal}(iii) implies that 
$$
(B-\theta_{\dbar-2}^*)(B-\theta_{\dbar-1}^*)A w_{0,\dbar-1}=0.
$$ 
Since $q^2\not=1$ and $w_{0,\dbar-1}$ is a basis for $W_\lambda^\delta(a,b,c)(\mu)$, there is a scalar $\varphi\in \F$ such that 
$
(B-\theta_{\dbar-2}^*)A w_{0,\dbar-1}=(q-q^{-1})\varphi\cdot w_{0,\dbar-1}.
$
The condition (ii) follows.

Applying Lemma \ref{lem:W_action}(i) a direct calculation shows that 
$$
\varphi =\frac{\varphi_{\dbar-1}+\theta_{\dbar-1}(\theta_{\dbar-1}^*-\theta_{\dbar-2}^*)}{q-q^{-1}}.
$$
By Definition \ref{defn:feasible} it is straightforward to verify that $(\bar{a},\bar{b},\bar{c},\bar{\lambda})$ is feasible for $(\mu,\varphi,\omega^*,\omega^\e)$ where 
$\omega^*$ and $\omega^{\e}$ are the scalars (\ref{omega*}) and (\ref{omegae}) respectively. By Theorem \ref{thm:feasible&universal} there exists a $\triangle_q$-module homomorphism 
$
M_{\bar{\lambda}}(\bar{a},\bar{b},\bar{c})\to W_\lambda^{\delta}(a,b,c)
$
that $\bar{m}_0$ to $w_{0,\dbar-1}$. 
Since $(\bar{a}^{\dbar}\bar{\lambda}^{-\dbar},\bar{\delta})=(a^{-\dbar}\lambda^{\dbar},\delta)$ the equation (\ref{delta&bardelta}) holds. 
By Proposition \ref{prop:delta&W} the condition (iii) follows. 
\end{proof}

\begin{prop}
\label{prop2:existence}
For any $(a,b,c,\lambda,\delta),(\bar{a},\bar{b},\bar{c},\bar{\lambda},\bar{\delta})\in {\F^\times}^4\times \F$ with $(a,b,c,\lambda,\delta)\sim (\bar{a},\bar{b},\bar{c},\bar{\lambda},\bar{\delta})$, there exists a $\triangle_q$-module homomorphism 
$W_{\bar{\lambda}}^{\bar{\delta}}(\bar{a},\bar{b},\bar{c})\to W_\lambda^\delta(a,b,c)$ that sends $\bar{w}_0$ to a marginal weight vector of $W_\lambda^\delta(a,b,c)$. 
\end{prop}
\begin{proof}
Immediate from Theorem \ref{thm:z2s4&W} and Lemmas \ref{lem1-0:existence} and \ref{lem1-2:existence}.
\end{proof}

\section{Proof for Theorem \ref{thm:bijection}}
\label{s:isoclass}

Recall the equivalence relation $\simeq$ on ${\F^\times}^4\times \F$ from Definition \ref{defn:simeq}. Recall the set ${\bf PM}_{\dbar}$ from below Definition \ref{defn:simeq}.

\begin{prop}
\label{prop:existence}
For any $(a,b,c,\lambda,\delta)\in {\bf PM}_{\dbar}$ and any $(\bar{a},\bar{b},\bar{c},\bar{\lambda},\bar{\delta})\in {\F^\times}^4\times \F$ with $(a,b,c,\lambda,\delta)\simeq (\bar{a},\bar{b},\bar{c},\bar{\lambda},\bar{\delta})$, the $\triangle_q$-module $W_\lambda^\delta(a,b,c)$ 
is isomorphic to 
$W_{\bar{\lambda}}^{\bar{\delta}}(\bar{a},\bar{b},\bar{c})$.
Moreover the set ${\bf PM}_{\dbar}$ is closed under $\simeq$.
\end{prop}
\begin{proof}
Immediate from Theorem \ref{thm2:irr} and Proposition \ref{prop2:existence}.
\end{proof}

\begin{lem}
\label{lem0:onto}
If $(a,b,c,\lambda), (\bar{a},\bar{b},\bar{c},\bar{\lambda})\in {\F^\times}^4$ with $(a,b,c,\lambda)\fz2s4\! (\bar{a},\bar{b},\bar{c},\bar{\lambda})$, then $\theta_i^*=\bar{\theta}_i^*$ for all $i\in \N$. 
\end{lem}
\begin{proof}
Since $(a,b,c,\lambda)\fz2s4\! (\bar{a},\bar{b},\bar{c},\bar{\lambda})$ and by Theorem \ref{thm:feasible}(ii), the elements $(a,b,c,\lambda)$ and $(\bar{a},\bar{b},\bar{c},\bar{\lambda})$ are feasible for the same element of $\F^\times \times \F^3$. 
Hence $b\lambda^{-1}=\bar{b}\bar{\lambda}^{-1}$ by Definition \ref{defn:feasible}(i). 
By (\ref{thetais}) the lemma follows.
\end{proof}

\begin{prop}
\label{prop0:onto}
Suppose that there is an integer $i$ with $1\leq i\leq \dbar-1$ such that $\varphi_i=0$. 
Then there exists an element $(\bar{a},\bar{b},\bar{c},\bar{\lambda},\bar{\delta})\in {\F^\times}^4\times \F$ with $(a,b,c,\lambda,\delta)\z2s4\! (\bar{a},\bar{b},\bar{c},\bar{\lambda},\bar{\delta})$ satisfying the following conditions:
\begin{enumerate}
\item $\bar{\lambda}^2=q^{2(i-1)}$.

\item $\bar{\theta}_h=\theta_h$ for all $h\in \N$.

\item $\bar{\theta}_h^*=\theta_h^*$ for all $h\in \N$.

\item $\bar{\varphi}_h=\varphi_h$ for all $h\in \N$.

\item There is a $\triangle_q$-module isomorphism $W_{\bar{\lambda}}^{\bar{\delta}}(\bar{a},\bar{b},\bar{c})\to W_\lambda^\delta(a,b,c)$ that maps 
\begin{eqnarray*}
\bar{w}_h &\mapsto &  w_h
\qquad 
\hbox{for all $h=0,1,\ldots,\dbar-1$}.
\end{eqnarray*}
\end{enumerate} 
\end{prop}
\begin{proof}
Since $\varphi_i=0$ it follows from (\ref{varphii}) that 
$q^{2(i-1)}\in\{\lambda^2, a^{-1}b^{-1}c^{-1}\lambda q^{-1},a^{-1}b^{-1}c\lambda q^{-1}\}$. 
We select $(\bar{a},\bar{b},\bar{c},\bar{\lambda})$ from the $2^{\rm nd}$, $4^{\rm th}$ and $7^{\rm th}$ rows of Table \ref{t:z2s4} as follows:
\begin{gather*}
(\bar{a},\bar{b},\bar{c},\bar{\lambda})
=\left\{
\begin{array}{ll}
(a,b,c,\lambda) 
\qquad &\hbox{if $q^{2(i-1)}=\lambda^2$},
 \\
(\frac{a}{\sqrt{abc\lambda q}},
\frac{b}{\sqrt{abc\lambda q}},
\frac{c}{\sqrt{abc\lambda q}},
\frac{\lambda}{\sqrt{abc\lambda q}}) 
\qquad &\hbox{if $q^{2(i-1)}=a^{-1}b^{-1}c^{-1}\lambda q^{-1}$},
 \\
(\frac{ac}{\sqrt{abc\lambda q}},
\frac{bc}{\sqrt{abc\lambda q}},
\frac{1}{\sqrt{abc\lambda q}},
\frac{c\lambda}{\sqrt{abc\lambda q}}) 
\qquad &\hbox{if $q^{2(i-1)}=a^{-1}b^{-1}c\lambda q^{-1}$}.
\end{array}
\right.
\end{gather*}
By construction the condition (i) holds. 
Since $\bar{a}\bar{\lambda}^{-1}=a\lambda^{-1}$ the condition (ii) holds by (\ref{thetai}). 
Since $(a,b,c,\lambda)\fz2s4\! (\bar{a},\bar{b},\bar{c},\bar{\lambda})$ by Definition \ref{defn:fz2s4&z2s4}(i), the condition (iii) is immediate from Lemma \ref{lem0:onto}.  
It is straightforward to verify the condition (iv) by using (\ref{varphii}).

Let $\bar{\delta}=\delta$. Since $(\bar{a}\bar{\lambda}^{-1},\bar{\delta})=(a\lambda^{-1},\delta)$ the equation (\ref{delta&bardelta}) holds. 
Hence $(a,b,c,\lambda,\delta)\z2s4\! (\bar{a},\bar{b},\bar{c},\bar{\lambda},\bar{\delta})$ by Definition \ref{defn:fz2s4&z2s4}(ii). 
By Theorem \ref{thm:z2s4&W} there exists a $\triangle_q$-module isomorphism 
$$
f:W_{\bar{\lambda}}^{\bar{\delta}}(\bar{a},\bar{b},\bar{c})\to W_\lambda^\delta(a,b,c)
$$ 
that maps $\bar{w}_0$ to $w_0$.  
Applying (ii) and Lemma \ref{lem:W_action}(i) yields that 
$f(\bar{w}_h)=w_h$ 
for all $h=0,1,\ldots,\dbar-1$. The condition (v) follows. 
\end{proof}

While using $(\bar{a},\bar{b},\bar{c},\bar{\lambda},\bar{\delta})$ to denote an element of ${\F^\times}^4\times \F$, the notation $\bar{w}_{ij}$ represents the vector of $W_{\bar{\lambda}}^{\bar{\delta}}(\bar{a},\bar{b},\bar{c})$ corresponding to the vector $w_{ij}$ of $W_\lambda^\delta(a,b,c)$ for any integers $i,j$ with $0\leq i,j\leq \dbar-1$.

\begin{prop}
\label{prop1:onto}
Suppose that there is an integer $i$ with $0\leq i\leq \dbar-1$ such that $\varphi_i=0$.  
Then there exists an element $(\bar{a},\bar{b},\bar{c},\bar{\lambda},\bar{\delta})\in {\F^{\times}}^4\times \F$ with $(a,b,c,\lambda,\delta)\simeq (\bar{a},\bar{b},\bar{c},\bar{\lambda},\bar{\delta})$ satisfying the following conditions:
\begin{enumerate}
\item $\bar{\theta}_h^*=\theta_{i+h}^*$ for all $h\in \N$.

\item There is a $\triangle_q$-module homomorphism $
f:W_{\bar{\lambda}}^{\bar{\delta}}(\bar{a},\bar{b},\bar{c})\to W_\lambda^\delta(a,b,c)$ 
that maps $\bar{w}_0$ to $w_i$.

\item 
$f(\bar{w}_{0,j-i})=w_{ij}$ for any integer $j$ with $i\leq j\leq \dbar-1$ and $\theta_j^*\not\in\{\theta_h^
*\,|\,h=i,i+1,\ldots,j-1\}$. 
\end{enumerate}
\end{prop}
\begin{proof}
If $i=0$ then the proposition follows by selecting $(\bar{a},\bar{b},\bar{c},\bar{\lambda},\bar{\delta})=(a,b,c,\lambda,\delta)$. Suppose that $i\geq 1$. 
By Proposition \ref{prop0:onto} we may assume that $\lambda^2=q^{2(i-1)}$ instead of $\varphi_i=0$. Combined with Lemma \ref{lem1-0:existence} the conditions (i) and (ii) follow.

Suppose that $j$ is an integer with $i\leq j\leq \dbar-1$ satisfying the hypothesis $\theta_j^*\not=\theta_h^
*$ for all $h=i,i+1,\ldots,j-1$.   
Applying Lemma \ref{lem:W_action}(i) the vector 
$f(\bar{w}_h)$ is equal to $w_{i+h}$ plus a linear combination of $w_i,w_{i+1},\ldots,w_{i+h-1}$ for all $h=0,1,\ldots,\dbar-i-1$. Hence $f(\bar{w}_{0,j-i})$ is a linear combination of $\{w_h\}_{h=i}^j$. In addition the coefficient of $w_j$ in $f(\bar{w}_{0,j-i})$ is equal to the coefficient of $\bar{w}_{j-i}$ in $\bar{w}_{0,j-i}$. 
By Lemma \ref{lem1:Wmarginal} and since $\bar{\theta}_{j-i}^*=\theta_j^*$ by (i) it follows that $(B-\theta_j^*)\bar{w}_{0,j-i}=0$. 
Applying $f$ to the above equation we have 
$(B-\theta_j^*)f(\bar{w}_{0,j-i})=0$.
By Lemma \ref{lem1-2:Wmarginal} there is a scalar $\e\in \F$ such that 
\begin{gather}
\label{prop1:onto-(iii)-f}
f(\bar{w}_{0,j-i})=\e w_{ij}.
\end{gather}
By (\ref{wij}) the coefficient of $w_j$ in $w_{ij}$ is equal to $\prod_{h=i}^{j-1}(\theta_j^*-\theta_h^*)\not=0$
and the coefficient of $\bar{w}_{j-i}$ in $\bar{w}_{0,j-i}$ is equal to 
$
\prod_{h=1}^{j-i}(\bar{\theta}_{j-i}^*-\bar{\theta}_{h-1}^*)$. 
By (i) both coefficients are equal. Comparing the coefficients of $w_j$ in both sides of (\ref{prop1:onto-(iii)-f}) yields that $\e=1$. The condition (iii) follows.
\end{proof}

\begin{lem}
\label{lem2-0:onto}
Suppose that there is an integer $i$ with $1\leq i\leq \dbar-1$ such that $\varphi_1\varphi_2\cdots\varphi_i\not=0$.
Then the following conditions are equivalent: 
\begin{enumerate}
\item $(B-\theta_{i+1}^*)(B-\theta_i^*)w_{0i}=0$.

\item $q^{2(i-1)}\in\{b^{-2}q^{-2}, 
ab^{-1}c\lambda q^{-1}, 
ab^{-1}c^{-1}\lambda q^{-1}\}$.
\end{enumerate}
\end{lem}
\begin{proof}
(ii) $\Rightarrow$ (i): Immediate from Lemma \ref{lem3:Wmarginal}.

(i) $\Rightarrow$ (ii): 
Looking at the coefficient of $w_0$ in $w_{0,i-1}$ by (\ref{wij}) yields that $w_{0,i-1}$ is nonzero.
Since $\varphi_i\not=0$ and by Lemma \ref{lem3:Wmarginal} the implication (i) $\Rightarrow$ (ii) follows. 
\end{proof}

Later the quintuple $(\hat{a},\hat{b},\hat{c},\hat{\lambda},\hat{\delta})$ will also be used to denote an element of ${\F^\times}^4\times \F$. We will adopt a similar convention for the notations
$\{\hat{\theta}_i\}_{i\in \N},\{\hat{\theta}_i^*\}_{i\in \N}, \{\hat{\varphi}_i\}_{i\in \N}, \{\hat{w}_i\}_{i=0}^{\dbar-1}, \{\hat{w}_{ij}\}_{0\leq i,j\leq \dbar-1}$.

\begin{prop}
\label{prop2:onto}
Suppose that there is an integer $i$ with $0\leq i\leq \dbar-1$ such that 
$\varphi_1\varphi_2\cdots\varphi_i\not=0$. If 
$$
(B-\theta_{i+1}^*)(B-\theta_i^*)w_{0i}=0,
$$  
then there exists an element $(\bar{a},\bar{b},\bar{c},\bar{\lambda},\bar{\delta})\in {\F^\times}^4\times \F$ with $(a,b,c,\lambda,\delta)\simeq (\bar{a},\bar{b},\bar{c},\bar{\lambda},\bar{\delta})$ satisfying the following conditions:
\begin{enumerate}
\item $\bar{\theta}_h^*=\theta_{i+h}^*$ 
for all $h\in \N$.

\item There is a $\triangle_q$-module homomorphism 
$
W_{\bar{\lambda}}^{\bar{\delta}}(\bar{a},\bar{b},\bar{c})\to W_\lambda^\delta(a,b,c)
$
that maps $\bar{w}_0$ to $w_{0i}$.
\end{enumerate}
\end{prop}
\begin{proof} 
If $i=0$ then the proposition follows by selecting $(\bar{a},\bar{b},\bar{c},\bar{\lambda},\bar{\delta})=(a,b,c,\lambda,\delta)$. 
Suppose that $i\geq 1$. 
It follows from Lemma \ref{lem2-0:onto} that $q^{2(i-1)}\in\{b^{-2}q^{-2}, 
ab^{-1}c\lambda q^{-1}, 
ab^{-1}c^{-1}\lambda q^{-1}\}$. 
We select $(\hat{a},\hat{b},\hat{c},\hat{\lambda})$ from the $10^{\rm th}$, $13^{\rm th}$ and $18^{\rm th}$ rows of Table \ref{t:z2s4} as follows:
\begin{gather*}
(\hat{a},\hat{b},\hat{c},\hat{\lambda})
=\left\{
\begin{array}{ll}
(c, \lambda^{-1}q^{-1}, a, b^{-1}q^{-1})
\qquad &\hbox{if $q^{2(i-1)}=b^{-2}q^{-2}$},
 \\
(\frac{c}{\sqrt{abc\lambda q}},
\frac{abc}{\sqrt{abc\lambda q}},
\frac{a}{\sqrt{abc\lambda q}},
\frac{ac\lambda}{\sqrt{abc\lambda q}}) 
\qquad &\hbox{if $q^{2(i-1)}=ab^{-1}c\lambda q^{-1}$},
 \\
(\frac{1}{\sqrt{abc\lambda q}},
\frac{ab}{\sqrt{abc\lambda q}},
\frac{ac}{\sqrt{abc\lambda q}},
\frac{a\lambda}{\sqrt{abc\lambda q}}) 
\qquad &\hbox{if $q^{2(i-1)}=ab^{-1}c^{-1}\lambda q^{-1}$}.
\end{array}
\right.
\end{gather*}
By construction $\hat{\lambda}^2=q^{2(i-1)}$. Note that $(a,b,c,\lambda)\fz2s4\! (\hat{a},\hat{b},\hat{c},\hat{\lambda})$ by Definition \ref{defn:fz2s4&z2s4}(i).
Set  
$$
\hat{\delta}=\delta+a^{\dbar}\lambda^{-\dbar}+a^{-\dbar}\lambda^{\dbar}-\hat{a}^{\dbar}\hat{\lambda}^{-\dbar}-\hat{a}^{-\dbar}\hat{\lambda}^{\dbar}.
$$
Then $(a,b,c,\lambda,\delta)\z2s4\! (\hat{a},\hat{b},\hat{c},\hat{\lambda},\hat{\delta})$ by Definition \ref{defn:fz2s4&z2s4}(ii). By Theorem \ref{thm:z2s4&W} there exists a $\triangle_q$-module isomorphism
$$
f_1:W_{\hat{\lambda}}^{\hat{\delta}}(\hat{a},\hat{b},\hat{c})\to 
W_\lambda^\delta(a,b,c)
$$
that maps $\hat{w}_0$ to $w_0$. Since $\hat{\lambda}^2=q^{2(i-1)}$ it follows from Lemma \ref{lem1-0:existence} that there exists an element $(\bar{a},\bar{b},\bar{c},\bar{\lambda},\bar{\delta})\in {\F^\times}^4\times \F$ with $(\hat{a},\hat{b},\hat{c},\hat{\lambda},\hat{\delta})\sim (\bar{a},\bar{b},\bar{c},\bar{\lambda},\bar{\delta})$ such that $\bar{\theta}_h^*=\hat{\theta}_{i+h}^*$ for all $h\in \N$ and there is a 
$\triangle_q$-module homomorphism 
$$
f_2:W_{\bar{\lambda}}^{\bar{\delta}}(\bar{a},\bar{b},\bar{c})\to 
W_{\hat{\lambda}}^{\hat{\delta}}(\hat{a},\hat{b},\hat{c})
$$
that maps $\bar{w}_0$ to $\hat{w}_i$. 
Since $(\hat{a},\hat{b},\hat{c},\hat{\lambda},\hat{\delta})\sim (\bar{a},\bar{b},\bar{c},\bar{\lambda},\bar{\delta})$ and $(a,b,c,\lambda,\delta)\sim (\hat{a},\hat{b},\hat{c},\hat{\lambda},\hat{\delta})$ by Definition \ref{defn:sim}(i) it follows from Definition \ref{defn:simeq} that 
$
(a,b,c,\lambda,\delta)\simeq (\bar{a},\bar{b},\bar{c},\bar{\lambda},\bar{\delta})$.

Since $\theta_h^*=\hat{\theta}_h^*$ for all $h\in \N$ by Lemma \ref{lem0:onto}, the condition (i) follows. To see (ii) it suffices to show that $(f_1\circ f_2)(\bar{w}_0)$ is a nonzero scalar multiple of $w_{0i}$. 
Since $\hat{\lambda}^2=q^{2(i-1)}$ it follows from (\ref{varphii}) that $\hat{\varphi}_i=0$.
Hence   
$$
(B-\theta_i^*)\hat{w}_i=0
$$
by Lemma \ref{lem1:Wmarginal}. 
Since $f_1(\hat{w}_0)=w_0$ and by Lemma \ref{lem:W_action}(i) the vector $f_1(\hat{w}_i)$ is a linear combination of $\{w_h\}_{h=0}^i$. 
By Lemma \ref{lem1-2:Wmarginal} there is a scalar $\e\in \F$ such that 
$$
(f_1\circ f_2)(w_0)=f_1(\hat{w}_i)=\e w_{0i}.
$$
Since $f_1$ is a $\triangle_q$-module isomorphism the vector $f_1(\hat{w}_i)$ is nonzero and the scalar $\e$ is nonzero. 
The proposition follows.
\end{proof}

\begin{lem}
\label{lem3-1:onto}
Suppose there are two integers $i,j$ with $0\leq i\leq j\leq \dbar-1$ such that $\varphi_i=0$ and $\theta_j^*\not\in\{\theta_h^
*\,|\,h=i,i+1,\ldots,j-1\}$. Then there exist an integer $k$ with $i\leq k\leq j$ an element $(\bar{a},\bar{b},\bar{c},\bar{\lambda},\bar{\delta})\in {\F^\times}^4\times \F$ with $(a,b,c,\lambda,\delta)\simeq (\bar{a},\bar{b},\bar{c},\bar{\lambda},\bar{\delta})$ satisfying the following conditions:
\begin{enumerate}
\item $\bar{\varphi}_1\bar{\varphi}_2\cdots\bar{\varphi}_{j-k}\not=0$.

\item $\bar{\theta}_h^*=\theta_{k+h}^*$ 
for all $h\in \N$.

\item There is a $\triangle_q$-module homomorphism $W_{\bar{\lambda}}^{\bar{\delta}}(\bar{a},\bar{b},\bar{c})\to W_\lambda^\delta(a,b,c)$ that maps $\bar{w}_{0,j-k}$ to $w_{ij}$.
\end{enumerate}
\end{lem}
\begin{proof}
We proceed by induction on $j-i$. If $i=j$ then the lemma is immediate from Proposition \ref{prop1:onto}. Now suppose that $j>i$. By Proposition \ref{prop1:onto} there exists an element $(\hat{a},\hat{b},\hat{c},\hat{\lambda},\hat{\delta})\in {\F^\times}^4\times \F$ with $(a,b,c,\lambda,\delta)\simeq (\hat{a},\hat{b},\hat{c},\hat{\lambda},\hat{\delta})$ satisfying the following conditions:
\begin{enumerate}
\item[(1)] $\hat{\theta}_h^*=\theta_{i+h}^*$ for all $h\in \N$.

\item[(2)] There is a $\triangle_q$-module homomorphism  $f:W_{\hat{\lambda}}^{\hat{\delta}}(\hat{a},\hat{b},\hat{c})\to W_\lambda^\delta(a,b,c)$ such that 
$$
f(\hat{w}_{0,j-i})=w_{ij}.
$$
\end{enumerate}
If $\hat{\varphi}_1\hat{\varphi}_2\cdots\hat{\varphi}_{j-i}\not=0$ then the lemma follows.
Thus suppose that there is an integer $\ell$ with $1\leq \ell\leq j-i$ such that $\hat{\varphi}_{\ell}=0$. By the induction hypothesis there exist an integer $m$ with $\ell\leq m\leq j-i$ and an element $(\bar{a},\bar{b},\bar{c},\bar{\lambda},\bar{\delta})\in {\F^\times}^4\times \F$ with $(\hat{a},\hat{b},\hat{c},\hat{\lambda},\hat{\delta})\simeq (\bar{a},\bar{b},\bar{c},\bar{\lambda},\bar{\delta})$ satisfying the following conditions:
\begin{enumerate}
\item[(a)] $\bar{\varphi}_1\bar{\varphi}_2\cdots\bar{\varphi}_{j-i-m}\not=0$.

\item[(b)] $\bar{\theta}_h^*=\hat{\theta}_{m+h}^*$ 
for all $h\in \N$.

\item[(c)] There is a $\triangle_q$-module homomorphism $g:W_{\bar{\lambda}}^{\bar{\delta}}(\bar{a},\bar{b},\bar{c})\to W_{\hat{\lambda}}^{\hat{\delta}}(\hat{a},\hat{b},\hat{c})$ such that 
$$
g(\bar{w}_{0,j-i-m})=\hat{w}_{\ell,j-i}.
$$
\end{enumerate}

By the transitive property of $\simeq$ 
the element $(a,b,c,\lambda,\delta)\simeq (\bar{a},\bar{b},\bar{c},\bar{\lambda},\bar{\delta})$. Set 
$$
k=i+m.
$$ 
Then the condition (i) is exactly the condition (a). The condition (ii) is immediate from the conditions (1) and (b).
Since $\hat{\varphi}_{\ell}=0$ and by (\ref{wij}) the vector $\hat{w}_{0,j-i}$ is a linear combination of $\{\hat{w}_h\}_{h=\ell}^{j-i}$. By Lemma \ref{lem1:Wmarginal},
$(B-\hat{\theta}_{j-i}^*)\hat{w}_{0,j-i}=0$. 
By Lemma \ref{lem1-2:Wmarginal} there is a scalar $\e\in \F$ such that $\hat{w}_{0,j-i}=\e \hat{w}_{\ell,j-i}$.
Applying the conditions (2) and (c) yields that $\e\cdot (f\circ g)$ satisfies the condition (iii). The lemma follows.
\end{proof}

\begin{prop}
\label{prop3:onto}
Suppose that $(a,b,c,\lambda,\delta)\in {\bf PM}_{\dbar}$ and there are two integers $i,j$ with $0\leq i\leq j\leq \dbar-1$ such that $\varphi_i=0$ and $\theta_j^*\not\in\{\theta_h^
*\,|\,h=i,i+1,\ldots,j-1\}$. If
\begin{gather}
\label{e3:onto}
(B-\theta_{j+1}^*)(B-\theta_j^*)w_{ij}=0,
\end{gather}
then there exists an element $(\bar{a},\bar{b},\bar{c},\bar{\lambda},\bar{\delta})\in {\F^\times}^4\times \F$ with $(a,b,c,\lambda,\delta)\simeq (\bar{a},\bar{b},\bar{c},\bar{\lambda},\bar{\delta})$ satisfying the following conditions:
\begin{enumerate}
\item $\bar{\theta}_h^*=\theta_{j+h}^*$ 
for all $h\in \N$.

\item There is a $\triangle_q$-module isomorphism 
$
W_{\bar{\lambda}}^{\bar{\delta}}(\bar{a},\bar{b},\bar{c})\to W_\lambda^\delta(a,b,c)
$
that maps $\bar{w}_0$ to $w_{ij}$.
\end{enumerate}
\end{prop}
\begin{proof}  
There are an integer $k$ with $i\leq k\leq j$ and an element $(\bar{a},\bar{b},\bar{c},\bar{\lambda},\bar{\delta})\in {\F^\times}^4\times \F$ with $(a,b,c,\lambda,\delta)\simeq (\bar{a},\bar{b},\bar{c},\bar{\lambda},\bar{\delta})$ 
satisfying Lemma \ref{lem3-1:onto}(i)--(iii).  
Let $f$ denote the $\triangle_q$-module homomorphism described in Lemma \ref{lem3-1:onto}(iii). Since $(a,b,c,\lambda,\delta)\in {\bf PM}_{\dbar}$ and by Theorem \ref{thm2:irr} the map $f$ is a $\triangle_q$-module isomorphism.  
Pulling back (\ref{e3:onto}) via $f$ yields that
$$
(B-\bar{\theta}_{j-k+1}^*)(B-\bar{\theta}_{j-k}^*)\bar{w}_{0,j-k}=0.
$$
Combined with Proposition \ref{prop2:onto} the proposition follows.
\end{proof}

\begin{prop}
\label{prop4:onto}
Suppose that there are two integers $i,j$ with $0\leq i<j\leq \dbar-1$ such that $\varphi_{i+1}=0$ and $\theta_i^*=\theta_j^*$. 
Then there exists an element $(\bar{a},\bar{b},\bar{c},\bar{\lambda},\bar{\delta})\in {\F^\times}^4\times \F$ with $(a,b,c,\lambda,\delta)\simeq (\bar{a},\bar{b},\bar{c},\bar{\lambda},\bar{\delta})$ satisfying the following conditions: 
\begin{enumerate}
\item $\bar{\theta}_h^*=\theta_{j+h}^*$ for all $h\in \N$.

\item If $j=i+1$ then the following equation holds:
$$
\frac{\varphi_i+\theta_i(\theta_i^*-\theta_{i+2}^*)}{q-q^{-1}}=(\bar{c}+\bar{c}^{-1})(\bar{\lambda}-\bar{\lambda}^{-1})-(\bar{a}+\bar{a}^{-1})(\bar{b}q-\bar{b}^{-1}q^{-1}).
$$
\end{enumerate} 
\end{prop}
\begin{proof}
By Proposition \ref{prop0:onto} we may assume that 
\begin{gather}
\label{prop4:onto-1}
\lambda^2=q^{2i}
\end{gather}
instead of $\varphi_{i+1}=0$. 
Since $\theta_i^*=\theta_j^*$ it follows that $\lambda^2=b^2q^{2(i+j)}$. Substituting (\ref{prop4:onto-1}) into the above equation yields that 
\begin{gather}
\label{prop4:onto-2}
b^{-2}q^{-2}=q^{2(j-1)}.
\end{gather}

Set 
$
(\bar{a},\bar{b},\bar{c},\bar{\lambda},\bar{\delta})
=
(a,b^{-1},c,\lambda,\delta).
$
Then  
$\bar{b}\bar{\lambda}^{-1}=b\lambda^{-1}q^{2j}$
by (\ref{prop4:onto-2}). The condition (i) holds.
Using (\ref{prop4:onto-1}) and (\ref{prop4:onto-2}) it is straightforward to verify the condition (ii). It remains to show that $(a,b,c,\lambda,\delta)\simeq (\bar{a},\bar{b},\bar{c},\bar{\lambda},\bar{\delta})$. Let 
\begin{align*}
\delta'=\delta+a^{\dbar}\lambda^{-\dbar}+a^{-\dbar}\lambda^{\dbar}-(b^{\dbar}c^{-\dbar}+b^{-\dbar}c^{\dbar})q^{\dbar}.
\end{align*}
By the $23^{\rm rd}$ row of Table \ref{t:z2s4} and  Definition \ref{defn:fz2s4&z2s4}(ii), 
$$
(a,b,c,\lambda,\delta)\z2s4 \!
(c^{-1},\lambda^{-1}q^{-1},a,b^{-1}q^{-1},\delta').
$$
Since (\ref{prop4:onto-2}) holds it follows from Definition \ref{defn:sim}(ii) that
$$
(c^{-1},\lambda^{-1}q^{-1},a,b^{-1}q^{-1},\delta')
\sim 
(c,\lambda^{-1}q^{-1},a,b q^{-1},\delta').
$$
By the $18^{\rm th}$ row of Table \ref{t:z2s4} and  Definition \ref{defn:fz2s4&z2s4}(ii), 
$$
(c,\lambda^{-1}q^{-1},a,b q^{-1},\delta')
\z2s4 \!
(\bar{a},\bar{b},\bar{c},\bar{\lambda},\bar{\delta}).
$$
By the transitive property of $\simeq$ the relation $(a,b,c,\lambda,\delta)\simeq (\bar{a},\bar{b},\bar{c},\bar{\lambda},\bar{\delta})$ holds. 
\end{proof}

Recall the statement of Theorem \ref{thm:bijection}.

\begin{proof}[Proof of Theorem \ref{thm:bijection}]
By Proposition \ref{prop:existence} the function from ${\bf PM}_{\dbar}/\!\simeq$ into ${\bf IM}_{\dbar}$ exists. By Theorem \ref{thm:=d} this function is onto. To see the injectivity, we assume that $(a,b,c,\lambda,\delta)$ and  $(\hat{a},\hat{b},\hat{c},\hat{\lambda},\hat{\delta})$ are in ${\bf PM}_{\dbar}$ and there exists a $\triangle_q$-module isomorphism 
$
f:W_{\hat{\lambda}}^{\hat{\delta}}(\hat{a},\hat{b},\hat{c})\to W_\lambda^\delta(a,b,c).
$ 
We need to show that 
\begin{gather}
\label{simeq}
(a,b,c,\lambda,\delta)\simeq (\hat{a},\hat{b},\hat{c},\hat{\lambda},\hat{\delta}).
\end{gather}  
We claim that there exists an element $(\bar{a},\bar{b},\bar{c},\bar{\lambda},\bar{\delta})\in {\F^\times}^4\times \F$ such that $(a,b,c,\lambda,\delta)\simeq 
(\bar{a},\bar{b},\bar{c},\bar{\lambda},\bar{\delta})$ and there exists a $\triangle_q$-module homomorphism
$$
g:W_{\bar{\lambda}}^{\bar{\delta}}(\bar{a},\bar{b},\bar{c})\to W_\lambda^\delta(a,b,c)
$$ 
that maps $\bar{w}_0$ to $f(\hat{w}_0)$. 
Suppose that the claim is true. Then the $\triangle_q$-module homomorphism 
$
f^{-1}\circ g:  
W_{\bar{\lambda}}^{\bar{\delta}}(\bar{a},\bar{b},\bar{c})
\to 
W_{\hat{\lambda}}^{\hat{\delta}}(\hat{a},\hat{b},\hat{c}) 
$
maps $\bar{w}_0$ to $\hat{w}_0$. Applying Theorem \ref{thm:z2s4&W} yields that $(\hat{a},\hat{b},\hat{c},\hat{\lambda},\hat{\delta})\z2s4\! (\bar{a},\bar{b},\bar{c},\bar{\lambda},\bar{\delta})$. In particular  
$(\hat{a},\hat{b},\hat{c},\hat{\lambda},\hat{\delta})\sim (\bar{a},\bar{b},\bar{c},\bar{\lambda},\bar{\delta})$ by Definition \ref{defn:sim}(i). 
Combined with the relation $(a,b,c,\lambda,\delta)\simeq 
(\bar{a},\bar{b},\bar{c},\bar{\lambda},\bar{\delta})$ the relation (\ref{simeq}) follows. Thus it suffices to prove the claim.

By Lemma \ref{lem:Verma_marginalweight}(i) the vector $\hat{w}_0$ is marginal weight vector of $W_{\hat{\lambda}}^{\hat{\delta}}(\hat{a},\hat{b},\hat{c})$ with weight 
$$
\mu=\hat{b}\hat{\lambda}^{-1}.
$$ 
Hence $f(\hat{w}_0)$ is a marginal weight vector of $W_\lambda^\delta(a,b,c)$ with weight $\mu$. 
Since any scalar appears in $\{\theta_i^*\}_{i=0}^{\dbar-1}$ at most twice,  
the dimension of $W_\lambda^\delta(a,b,c)(\mu)$ is less than or equal to two. 
Let 
$$
K=\ker (B-\mu q^2-\mu^{-1} q^{-2})(B-\mu-\mu^{-1}) A|_{W_\lambda^\delta(a,b,c)(\mu)}.
$$ 
Since $f(\hat{w}_0)\in K$ the dimension of $K$ is at least one. 
By Lemma \ref{lem:BeigenW} there is an integer $i$ with $0\leq i\leq \dbar-1$ such that $\theta_i^*=\mu+\mu^{-1}$. Hence $\mu\in\{b\lambda^{-1}q^{2i},b^{-1}\lambda q^{-2i}\}$ and 
\begin{gather*}
\mu q^2+\mu^{-1} q^{-2}=
\left\{
\begin{array}{ll}
\theta_{i+1}^* 
\qquad &\hbox{if $\mu=b\lambda^{-1}q^{2i}$},
\\
\theta_{i-1}^*
\qquad &\hbox{if $\mu=b^{-1}\lambda q^{-2i}$ and $i\not=0$},
\\
\theta_{\dbar-1}^*
\qquad &\hbox{if $\mu=b^{-1}\lambda q^{-2i}$ and $i=0$}.
\end{array}
\right.
\end{gather*}
To prove the claim, we divide the argument into the four cases:
(a) $\dim W_\lambda^\delta(a,b,c)(\mu)=\dim K=1$ and $\mu=b\lambda^{-1} q^{2i}$; 
(b) $\dim W_\lambda^\delta(a,b,c)(\mu)=\dim K=1$ and $\mu=b^{-1}\lambda q^{-2i}$; 
(c) $\dim W_\lambda^\delta(a,b,c)(\mu)=2>\dim K=1$; 
(d) $\dim W_\lambda^\delta(a,b,c)(\mu)=\dim K=2$.

(a): 
Applying Proposition \ref{prop1:onto} repeatedly we may assume that $\varphi_1\varphi_2\cdots\varphi_i\not=0$. 
By (\ref{wij}) the coefficient of $w_0$ in $w_{0i}$ is nonzero. Hence $w_{0i}$ is nonzero. Combined with Lemma \ref{lem1:Wmarginal} the vector $w_{0i}$ is a basis for
$W_\lambda^\delta(a,b,c)(\mu)$. Since $K=W_\lambda^\delta(a,b,c)(\mu)$ the vector $f(\hat{w}_0)$ is a nonzero scalar multiple of $w_{0i}$. The claim is now immediate from Proposition \ref{prop2:onto}.

(b): Since the former case has been done,  there is nothing to prove whenever $\mu$ is also equal to $b\lambda^{-1} q^{2j}$ for some integer $j$ with $0\leq j\leq \dbar-1$.
Thus we suppose that $\mu\not\in\{b\lambda^{-1} q^{2j}\,|\,j=0,1,\ldots,\dbar-1\}$. In other words 
\begin{gather}
\label{caseb}
\lambda^2\not=b^2 q^{2(i+j)}
\qquad 
\hbox{for all $j=0,1,\ldots,\dbar-1$}.
\end{gather}
Then $\theta_i^*\not=\theta_h^*$ for all $h=0,1,\ldots,i-1$. By (\ref{wij}) the coefficient of $w_i$ in $w_{0i}$ is nonzero. Combined with Lemma \ref{lem1:Wmarginal} the vector $w_{0i}$ is a basis for $W_\lambda^\delta(a,b,c)(\mu)$ as well as $K$. 

Suppose that $0\leq i\leq \dbar-2$. 
Then $\theta_{i+1}^*\not=\theta_h^*$ for all $h=0,1,\ldots,i$ by (\ref{caseb}).
By (\ref{wij}) the coefficient of $w_{i+1}$ in $w_{0,i+1}$ is nonzero. Hence $w_{0,i+1}$ is nonzero. 
Since $w_{0i}\in K$ it follows from Lemma \ref{lem2:Wmarginal}(i), (ii) that 
\begin{align*}
\left\{
\begin{array}{ll}
\theta_0^*=\theta_i^*
\qquad &\hbox{if $i\not=0$},
\\
\theta_1^*=\theta_{\dbar-1}^*
\qquad &\hbox{if $i=0$},
\end{array}
\right.
\end{align*}
which contradicts (\ref{caseb}). Therefore $i=\dbar-1$. 
Then  the condition (a) of Definition \ref{defn:sim}(iii) follows from (\ref{caseb}).
Since $w_{0,\dbar-1}\in K$ and $\theta_0^*\not=\theta_{\dbar-1}^*$ by (\ref{caseb}), it follows from Lemma \ref{lem2:Wmarginal}(iii) the condition (b) of Definition \ref{defn:sim}(iii) follows.  
Since $f(\hat{w}_0)$ is a nonzero scalar multiple of $w_{0,\dbar-1}$,  
the claim is immediate from Lemma \ref{lem1-2:existence}.

Before launching into the cases (c) and (d), we have some comments on 
$$
\dim W_{\lambda}^\delta(a,b,c)(\mu)=2.
$$ 
In this case there are two integers $i$ and $j$ with $0\leq i<j\leq \dbar-1$ such that $\theta_i^*=\theta_j^*=\mu+\mu^{-1}$. Then 
$\theta_h^*=\theta_{\ell}^*$
if and only if
$q^{2(h+\ell)}=q^{2(i+j)}$
for all distinct $h,\ell\in\{0,1,\ldots, \dbar-1\}$. It follows that 
\begin{enumerate}
\item[(1)] $\theta_i^*\not=\theta_h^*$ for all $h\in\{0,1,\ldots,\dbar-1\}\setminus\{i,j\}$;

\item[(2)] $\theta_j^*\not=\theta_h^*$ for all $h\in\{0,1,\ldots,\dbar-1\}\setminus\{i,j\}$;

\item[(3)] $\theta_{i+1}^*\not=\theta_h^*$ for all $h=0,1,\ldots,i-1$;

\item[(4)] $\theta_1^*\not=\theta_{\dbar-1}^*$ provided that $i=0$.
\end{enumerate} 
By the rank-nullity theorem there is an integer $k$ with $i<k\leq j$ such that $\varphi_k=0$. 
By (1) the coefficient of $w_i$ in $w_{0i}$ is nonzero. 
By (2) the coefficient of $w_j$ in $w_{kj}$ is nonzero. 
Then $w_{0i}$ and $w_{kj}$ are linearly independent. Combined with Lemma \ref{lem1:Wmarginal} the vectors $w_{0i}$ and $w_{kj}$ give a basis for $W_{\lambda}^\delta(a,b,c)(\mu)$.

Observe that  
$
\mu=b\lambda^{-1}q^{2i}=b^{-1}\lambda q^{-2j}$ or 
$
\mu=b^{-1}\lambda q^{-2i}=b\lambda^{-1}q^{2j}$. 
Applying Proposition \ref{prop1:onto} to $\varphi_k=0$ if necessary, we may only assume the latter case. Then
$$
\mu q^2+\mu^{-1} q^{-2}=\theta_{j+1}^* 
=
\left\{
\begin{array}{ll}
\theta_{i-1}^*
\qquad &\hbox{if $i\not=0$},
\\
\theta_{\dbar-1}^*
\qquad &\hbox{if $i=0$}.
\end{array}
\right.
$$
Furthermore the following conditions are equivalent:
%\begin{multicols}{2}
\begin{enumerate}
\item[(i)] $\dim K=2$.

\item[(ii)] $w_{0i}\in  K$.

\item[(iii)] $w_{0,i+1}=0$.

\item[(iv)] $\theta_i^*=\theta_{i+1}^*$. 

\item[(v)] $j=k=i+1$.
\end{enumerate}
%\end{multicols}
The proof of the equivalence (i)--(v) is as follows:
The implication (i) $\Rightarrow$ (ii) is obvious. 
In view of (1) and (4), Lemma \ref{lem2:Wmarginal}(i), (ii) implies the equivalence of (ii) and (iii). 
By (3) and (\ref{wij}) the equivalence of (iii)--(v) follows. 
Suppose that (v) holds. Then $\varphi_{i+1}=0$ and $w_{kj}=w_{i+1}$. It follows that $w_{i+1}\in K$. The implication (ii), (v) $\Rightarrow$ (i) follows.

(c): Since (ii) fails there is a scalar $\e\in \F$ such that $\e w_{0i}+w_{kj}$ is a basis for $K$. Then $f(\hat{w}_0)$ is a nonzero scalar multiple of $\e w_{0i}+w_{kj}$.
If $\e=0$ then the claim is immediate from Proposition \ref{prop3:onto}.

Suppose that $\e\not=0$. 
Since (ii) fails it follows from Lemma \ref{lem2:Wmarginal}(i), (ii) that $(B-\mu q^2-\mu^{-1}q^{-2})(B-\mu-\mu^{-1})A w_{0i}$ is a nonzero scalar multiple of $w_{0,i+1}$. Since (iv) fails and by (3) the coefficient of $w_{i+1}$ in $w_{0,i+1}$ is nonzero.
Since $(B-\mu q^2-\mu^{-1}q^{-2})(B-\mu-\mu^{-1})A w_{kj}$ a linear combination of $\{w_h\}_{h=k}^{j-1}$, this forces that $k=i+1$. 
Then there exists an element $(\bar{a},\bar{b},\bar{c},\bar{\lambda},\bar{\delta})\in {\F^\times}^4\times \F$ with $(a,b,c,\lambda,\delta)\simeq (\bar{a},\bar{b},\bar{c},\bar{\lambda},\bar{\delta})$ satisfying Proposition \ref{prop4:onto}(i).
By Proposition \ref{prop:existence} there exists a $\triangle_q$-module isomorphism 
$$
g:W_{\bar{\lambda}}^{\bar{\delta}}(\bar{a},\bar{b},\bar{c})\to W_\lambda^\delta(a,b,c).
$$
Since $\dim K=1$ the vector $g(\bar{w}_0)$ must be a nonzero scalar multiple of $f(\hat{w}_0)$. The claim follows.

(d): By (v) the scalar $\varphi_{i+1}=0$ and $w_{kj}=w_{i+1}$. 
Since $K=W_\lambda^\delta(a,b,c)(\mu)$ the space $W_\lambda^\delta(a,b,c)(\mu)$ is $(B-\mu q^2-\mu^{-1}q^{-2})A$-invariant. 
Set  
$$
\varphi=\frac{\varphi_i+\theta_i(\theta_i^*-\theta_{i+2}^*)}{q-q^{-1}}.
$$
Using Lemma \ref{lem:W_action}(i) a direct calculation shows that 
the $2\times 2$ matrix representing $(B-\mu q^2-\mu^{-1}q^{-2})A$ with respect to the ordering basis $w_{i+1}, w_{0i}$ for $W_\lambda^\delta(a,b,c)(\mu)$ is 
$$
\begin{pmatrix}
\varphi_{i+2}+\theta_{i+1}(\theta_i^*-\theta_{i+2}^*)

&(\theta_i^*-\theta_{i+2}^*)\prod\limits_{h=0}^{i-1}(\theta_i^*-\theta_h^*) 
\\
0 &(q-q^{-1})\varphi
\end{pmatrix}.
$$

The eigenvalues of $(B-\mu q^2-\mu^{-1}q^{-2})A$ on $W_\lambda^\delta(a,b,c)(\mu)$ are 
\begin{gather}
\varphi_{i+2}+\theta_{i+1}(\theta_i^*-\theta_{i+2}^*),
\label{K=2:eigen1}
\\
(q-q^{-1})\varphi.
\label{K=2:eigen2}
\end{gather}
By (1) the upper right entry of the $2\times 2$ matrix is nonzero.
Hence the geometric multiplicities of (\ref{K=2:eigen1}) and (\ref{K=2:eigen2}) are equal to one. 
By Definition \ref{defn:marginalvector} the vector $f(\hat{w}_0)$ is an eigenvector of $(B-\mu q^2-\mu^{-1}q^{-2})A$ in $W_\lambda^\delta(a,b,c)(\mu)$. 
If  the eigenvalue corresponding to $f(\hat{w}_0)$ is (\ref{K=2:eigen1}), then $f(\hat{w}_0)$ is a nonzero scalar multiple of $w_{i+1}$ and the claim follows from Proposition \ref{prop1:onto}. 
If the eigenvalue corresponding to $f(\hat{w}_0)$ is (\ref{K=2:eigen2}),  then the claim follows from Propositions \ref{prop:existence} and \ref{prop4:onto}.

We have shown that the claim is true in all cases. Theorem \ref{thm:bijection} is established.
\end{proof}

We finish the paper with a question.

\begin{prob}
Consider the Askey--Wilson algebras or their central extensions corresponding to 
hypergeometric orthogonal polynomials. 
For instance the Krawtchouk algebras, the universal Hahn algebra and the universal Racah algebra. 
Assume the underlying field is algebraically closed of positive characteristic. 
Please define what the marginal weights of their modules are and classify their finite-dimensional irreducible modules with marginal weights up to isomorphism.
\end{prob}

\newpage

\appendix

\section{The right $\S_4$-action of $\{\pm 1\}\backslash {\F^\times}^4$}
\label{s:S4}

Recall the left $\{\pm 1\}$-action on ${\F^\times}^4$ and the right $\S_4$-action on $\{\pm 1\}\backslash {\F^\times}^4$ from above Definition \ref{defn:fz2s4&z2s4}. 
Let $(a,b,c,\lambda)\in {\F^\times}^4$ be given. 
In the table below, we list every element $\sigma\in \S_4$ and a corresponding element $(\bar{a},\bar{b},\bar{c},\bar{\lambda})\in {\F^\times}^4$ satisfying 
$$
(\{\pm 1\}\cdot (a,b,c,\lambda))\cdot \sigma
=
\{\pm 1\}\cdot (\bar{a},\bar{b},\bar{c},\bar{\lambda}).
$$

\begin{table}[H]
\centering
\extrarowheight=1.5pt
\renewcommand{\arraystretch}{1.2}
\begin{tabular}{c||c|c|c}
$\sigma$  
&$(\bar{a},\bar{b},\bar{c},\bar{\lambda})$
&$\bar{\lambda}^2$ 
&$\bar{a}\bar{\lambda}^{-1}$
\\

\hline 
\hline

$1$
&$(a,b,c,\lambda)$ 
& \multirow{2}{*}{$\lambda^2$}
& \multirow{6}{*}{$a\lambda^{-1}$}
\\
\cline{1-2}
$(1\,2)$
&$(a,b,c^{-1},\lambda)$ 
& 
& 
\\
\cline{1-3}
$(2\,3)$
&$(\frac{a}{\sqrt{abc\lambda q}},
\frac{b}{\sqrt{abc\lambda q}},
\frac{c}{\sqrt{abc\lambda q}},
\frac{\lambda}{\sqrt{abc\lambda q}})$ 
& \multirow{2}{*}%{$a^{-1}b^{-1}c^{-1}\lambda q^{-1}$}  
{$\displaystyle \frac{\lambda}{abcq}$}
&
\\
\cline{1-2}
$(1\,3\,2)$
&$(\frac{a}{\sqrt{abc\lambda q}},
\frac{b}{\sqrt{abc\lambda q}},
\frac{\sqrt{abc\lambda q}}{c},
\frac{\lambda}{\sqrt{abc\lambda q}})$ 
& 
& 
\\
\cline{1-3}
$(1\,3)$
&$
(\frac{ac}{\sqrt{abc\lambda q}},
\frac{bc}{\sqrt{abc\lambda q}},
\sqrt{abc\lambda q},
\frac{c\lambda}{\sqrt{abc\lambda q}})$ 
& \multirow{2}{*}%{$a^{-1}b^{-1}c\lambda q^{-1}$}
{$\displaystyle \frac{c\lambda}{abq}$}
&  
\\
\cline{1-2}
$(1\,2\,3)$
&$
(\frac{ac}{\sqrt{abc\lambda q}},
\frac{bc}{\sqrt{abc\lambda q}},
\frac{1}{\sqrt{abc\lambda q}},
\frac{c\lambda}{\sqrt{abc\lambda q}})$ 
& 
& 
\\
 \midrule[1.5pt]
$(3\,4)$
& $(a^{-1},b,c,\lambda)$
& \multirow{2}{*}{$\lambda^2$}
& \multirow{6}{*}{$a^{-1}\lambda^{-1}$}
\\
\cline{1-2}
$(1\,2)(3\,4)$
&$(a^{-1},b,c^{-1},\lambda)$
& 
&
\\
\cline{1-3}
$(2\,4\,3)$
&$(
\frac{1}{\sqrt{a bc\lambda q}},
\frac{a b}{\sqrt{abc\lambda q}},
\frac{a c}{\sqrt{abc\lambda q}},
\frac{a \lambda}{\sqrt{a bc\lambda q}})$
&\multirow{2}{*}%{$ab^{-1}c^{-1}\lambda q^{-1}$}
{$\displaystyle \frac{a\lambda}{bcq}$}
&
\\
\cline{1-2}
$(1\,4\,3\,2)$
&$(
\frac{1}{\sqrt{a bc\lambda q}},
\frac{a b}{\sqrt{abc\lambda q}},
\frac{\sqrt{abc\lambda q}}{a c},
\frac{a \lambda}{\sqrt{a bc\lambda q}})$
&
&
\\
\cline{1-3}
$(1\,4\,3)$
&$(
\frac{c}{\sqrt{a bc\lambda q}},
\frac{a b c}{\sqrt{abc\lambda q}},
\frac{\sqrt{abc\lambda q}}{a},
\frac{a c\lambda}{\sqrt{a bc\lambda q}})$
&\multirow{2}{*}%{$ab^{-1}c\lambda q^{-1}$}
{$\displaystyle \frac{ac\lambda}{bq}$}
&
\\
\cline{1-2}
$(1\,2\,4\,3)$
&$(
\frac{c}{\sqrt{a bc\lambda q}},
\frac{a b c}{\sqrt{abc\lambda q}},
\frac{a}{\sqrt{abc\lambda q}},
\frac{a c\lambda}{\sqrt{a bc\lambda q}})$
& 
&
\\
 \midrule[1.5pt]
$(2\,4)$
&$(\sqrt{abc\lambda q},
\frac{ab}{\sqrt{abc\lambda q}},
\frac{ac}{\sqrt{abc\lambda q}},
\frac{a\lambda}{\sqrt{abc\lambda q}})$
& \multirow{2}{*}%{$ab^{-1}c^{-1}\lambda q^{-1}$}
{$\displaystyle \frac{a\lambda}{bcq}$}
&  \multirow{6}{*}{$bcq$}
\\
\cline{1-2}
$(1\,4\,2)$
&$(\sqrt{abc\lambda q},
\frac{ab}{\sqrt{abc\lambda q}},
\frac{\sqrt{abc\lambda q}}{ac},
\frac{a\lambda}{\sqrt{abc\lambda q}})$
& 
&
\\
\cline{1-3}
$(2\,3\,4)$
&$
(\frac{\sqrt{abc\lambda q}}{a},
\frac{b}{\sqrt{abc\lambda q}},
\frac{c}{\sqrt{abc\lambda q}},
\frac{\lambda}{\sqrt{abc\lambda q}})$
&\multirow{2}{*}%{$a^{-1}b^{-1}c^{-1}\lambda q^{-1}$}
{$\displaystyle \frac{\lambda}{abcq}$}
&
\\
\cline{1-2}
$(1\,3\,4\,2)$
&$
(\frac{\sqrt{abc\lambda q}}{a},
\frac{b}{\sqrt{abc\lambda q}},
\frac{\sqrt{abc\lambda q}}{c},
\frac{\lambda}{\sqrt{abc\lambda q}})$
&
&
\\
\cline{1-3}
$(1\,3)(2\,4)$
&$(c,\lambda^{-1}q^{-1},a,b^{-1}q^{-1})$
& \multirow{2}{*}%{$b^{-2}q^{-2}$}
{$\displaystyle \frac{1}{b^2q^2}$}
&
\\
\cline{1-2}
$(1\,4\,2\,3)$
&$(c,\lambda^{-1}q^{-1},a^{-1},b^{-1}q^{-1})$
& 
&

\\
\midrule[1.5pt]
$(1\,4)$
&$(\frac{\sqrt{abc\lambda q}}{c},
\frac{abc}{\sqrt{abc\lambda q}},
\frac{\sqrt{abc\lambda q}}{a},
\frac{ac\lambda}{\sqrt{abc\lambda q}})$
& \multirow{2}{*}%{$ab^{-1}c\lambda q^{-1}$}
{$\displaystyle \frac{ac\lambda}{bq}$}
&\multirow{6}{*}{$bc^{-1}q$}
\\
\cline{1-2}
$(1\,2\,4)$
&$(\frac{\sqrt{abc\lambda q}}{c},
\frac{abc}{\sqrt{abc\lambda q}},
\frac{a}{\sqrt{abc\lambda q}},
\frac{ac\lambda}{\sqrt{abc\lambda q}})$
& 
&
\\
\cline{1-3}
$(1\,4)(2\,3)$
&$(c^{-1},
\lambda^{-1}q^{-1},
a^{-1},
b^{-1}q^{-1})$
&\multirow{2}{*}%{$b^{-2}q^{-2}$} 
{$\displaystyle \frac{1}{b^2q^2}$}
&
\\
\cline{1-2}
$(1\,3\,2\,4)$
&$(c^{-1},
\lambda^{-1}q^{-1},
a,
b^{-1}q^{-1})$
& 
&
\\
\cline{1-3}
$(1\,3\,4)$
&$(\frac{\sqrt{abc\lambda q}}{ac},
\frac{bc}{\sqrt{abc\lambda q}},
\sqrt{abc\lambda q},
\frac{c\lambda}{\sqrt{abc\lambda q}})$
&\multirow{2}{*}
{$\displaystyle \frac{c\lambda}{abq}$}
&
\\
\cline{1-2}
$(1\,2\,3\,4)$
&$(\frac{\sqrt{abc\lambda q}}{ac},
\frac{bc}{\sqrt{abc\lambda q}},
\frac{1}{\sqrt{abc\lambda q}},
\frac{c\lambda}{\sqrt{abc\lambda q}})$
& 
&
\end{tabular}
\caption{The $\S_4$-orbit of $\{\pm 1\}\cdot (a,b,c,\lambda)$}\label{t:z2s4}
\end{table}

\subsection*{Acknowledgements}
The research was supported by the National Science and Technology Council of Taiwan under the projects MOST 105-2115-M-008-013, MOST 106-2628-M-008-001-MY4, MOST 110-2115-M-008-008-MY2 and NSTC 112-2115-M-008-009-MY2.

\bibliographystyle{amsplain}
\bibliography{MP}

\end{document}